\documentclass[leqno,12pt]{article}  
\usepackage{amsmath}
\usepackage{amssymb}

\usepackage[german,english]{babel}
\usepackage{amsthm}

\title{Compactifications of Log Morphisms}

\author{\textsc{Elmar Grosse-Kl\"onne}}
\date{}

\theoremstyle{plain} 
\newtheorem{satz}{Theorem}[section]  
\newtheorem{lem}[satz]{Lemma}  
\newtheorem{pro}[satz]{Proposition}  
 
\newcommand{\spec}{\mbox{\rm Spec}}  
\newcommand{\proj}{\mbox{\rm Proj}}
\newcommand{\quot}{\mbox{\rm Quot}}  
\newcommand{\spf}{\mbox{\rm Spf}}  

\newcommand{\bi}{\mbox{\rm Im}}  
\newcommand{\ke}{\mbox{\rm Ker}}  
\newcommand{\koke}{\mbox{\rm Coker}}  

\newcommand{\kara}{\mbox{\rm char}}  

\newcommand{\sym}{\mbox{\rm Sym}}

\newcommand{\dlog}{\mbox{\rm dlog}}
\newcommand{\Gr}{\mbox{\rm Gr}}

\theoremstyle{remark}

\theoremstyle{definition}

\begin{document}
\maketitle
\footnote[0]{Part of this work was done during my visit at the University of California, Berkeley, supported by the Deutsche Forschungs Gemeinschaft. I thank Arthur Ogus and Martin Olsson for helpful discussions and remarks. I am also very grateful to the referee for his careful reading and helpful suggestions which significantly improved the present paper.}

\begin{abstract}
We introduce the notion of a relative log scheme with boundary: a morphism of log schemes together with a (log schematically) dense open immersion of its source into a third log scheme. The sheaf of relative log differentials naturally extends to this compactification and there is a notion of smoothness for such data. We indicate how this weak sort of compactification may be used to develop useful de Rham and crystalline cohomology theories for semistable log schemes over the log point over a field which are not necessarily proper.    
\end{abstract}

%


\begin{center} {\bf Introduction} 
\end{center}

Let $X$ be a smooth variety over a field $k$. It is well known that for the study of the cohomology of $X$ --- or even for its very definition (e.g. crystalline \cite{maghreb}, rigid \cite{berco}), or the definition of nice coefficients for it (e.g. integrable connections with regular singularities) --- it is often indispensable to take into account also a boundary $D=\overline{X}-X$ of $X$ in a smooth compactification $X\subset\overline{X}$ of $X$. If $D\subset\overline{X}$ is a normal crossing divisor on $\overline{X}$, the cohomology can conveniently be studied in the framework of logarithmic algebraic geometry. On the other hand, log geometry proved also useful to define the cohomology of proper normal crossing varieties $X$ over $k$ which occur as a fibre of a semistable family, or more generally are $d$-semistable (\cite{fkato}), see \cite{steen}, \cite{mokr}. In the present paper we attempt to develop a concept in log geometry particularly suitable to treat the mixed situation: given a non-proper $d$-semistable normal crossing variety $X/k$, we want to explain how an open immersion of $X$ into a proper $k$-scheme $\overline{X}$ can be used to investigate the cohomology of $X$, the stress lying on the fact that we avoid the assumption that $\overline{X}$ be $d$-semistable and require a weaker condition instead.

Fix a base scheme $W$ for all occuring schemes. Let $T$ be a log scheme. The central definition of this note is that of a {\it $T$-log scheme with boundary}: A morphism of log schemes $X\to T$ together with an open log schematically dense embedding of log schemes $i:X\to\overline{X}$. For brevity, we often denote it simply by $(\overline{X},X)$. Morphisms of $T$-log schemes with boundary are defined in an obvious way. There are notions of exact and of boundary exact closed immersions of $T$-log schemes with boundary. The relative logarithmic de Rham complex $\Omega^{\bullet}_{X/T}$ on $X$ extends canonically to a complex $\Omega^{\bullet}_{(\overline{X},X)/T}$ on $\overline{X}$. These definitions are justified by a theory of smoothness for $T$-log schemes with boundary, well suited for cohomology purposes. Roughly, a $T$-log scheme with boundary $(\overline{X},X)$ is said to be weakly smooth if it satisfies a lifting property for morphisms from a nilpotent exact closed immersion of $T$-log schemes with boundary to $(\overline{X},X)$. Weak smoothness implies that $\Omega^{\bullet}_{(\overline{X},X)/T}$ is locally free. $(\overline{X},X)$ is said to be smooth if it is weakly smooth and if for boundary exact closed immersions $(\overline{Y},Y)\to(\overline{V},V)$ of $T$-log schemes with boundary, and morphisms $(\overline{Y},Y)\to(\overline{X},X)$ of $T$-log schemes with boundary, the projections $\overline{X}\overline{\times}_T\overline{V}\to\overline{V}$ lift log \'{e}tale locally near the image of $\overline{Y}$ in $\overline{X}\overline{\times}_T\overline{V}$ to strict and smooth morphisms of log schemes (see the text for the definition of $\overline{X}\overline{\times}_T\overline{V}$). This definition is of course geared to its application to (crystalline) cohomology. However, our main theorem gives a convenient criterion for smoothness in terms of morphisms of monoids, very similar to Kato's criterion for usual log smoothness. We emphasize that even if $f:X\to T$ actually extends to a morphism of log schemes $\overline{f}:\overline{X}\to T$, our notion of smoothness is more general: $(\overline{X},X)$ might be smooth as a $T$-log scheme with boundary while $\overline{f}$ is not a log smooth morphism in the usual sense (or even not ideally smooth as defined by Ogus \cite{ogid}). See for example the discussion at the beginning of Section 3. In this regard, the theme of this paper is that (usual) log smoothness in an `interior' $X\subset\overline{X}$ of a morphism of log schemes $\overline{f}:\overline{X}\to T$ should already ensure that $\overline{f}$ has nice cohomology. (A similar principle underlies the definition of rigid cohomology \cite{berco}.) We hope that our definitions are useful for a definition of log rigid cohomology, in the case of nontrivial log structures on the base; in special cases they already turned out to be so, see \cite{hkstrat}.

Section 1 contains the basic definitions and presents several examples. The main Section is the second one which is devoted to smoothness. The main theorem is the smoothness criterion \ref{tglatt}. In Section 3 we discuss the example of semistable $k$-log schemes with boundary (here $T$ is the log point over a field). These are smooth in the sense of Section 2 and we try to demonstrate how they can be used as substitutes for compactifications by usual semistable proper $k$-log schemes. We indicate several applications to de Rham cohomology and crystalline cohomology.\\

\section{$T$-log schemes with boundary}

\addtocounter{satz}{1}{\bf \arabic{section}.\arabic{satz}}\newcounter{basnot1}\newcounter{basnot2}\setcounter{basnot1}{\value{section}}\setcounter{basnot2}{\value{satz}} We fix a base scheme $W$; all schemes and morphisms of schemes are to be understood over $W$. All morphisms of schemes are quasi-separated. We also assume that all morphisms of schemes are quasi-compact: the only reason for this additional assumption is that it implies the existence of schematic images (=``closed images'') of morphisms: see \cite{EGA} I, 9.5. We say that an open immmersion $i:X\to\overline{X}$ is schematically dense if $\overline{X}$ coincides with the schematic image of $i$. For the basic notions of log algebraic geometry we refer to K. Kato \cite{kalo}. Log structures are understood for the \'{e}tale topology. By abuse of notation, for a scheme $X$ and a morphism of monoids $\alpha:N\to{\cal O}_X(X)$ (where ${\cal O}_X(X)$ is understood multiplicatively), we will denote by $(X,\alpha)$ the log scheme with underlying scheme $X$ whose log structure is associated with the chart $\alpha$. For a log scheme $(X,{\cal N}_X)=(X,{\cal N}_X\to{\cal O}_X)$ we will often just write $X$ if it is clear from the context to which log structure on $X$ we refer, i.e. in those cases the log structure is dropped in our notation. Similarly for morphisms of log schemes. An {\it exactification} of a closed immersion of fine log schemes $Y\to X$ is a factorization $Y\stackrel{i}{\to}Z\stackrel{f}{\to}X$ with $i$ an exact closed immersion and $f$ log \'{e}tale. Recall that a morphism of log schemes $f:(X,{\cal N}_X)\to(Y,{\cal N}_Y)$ is said to be {\it strict} if $f^*{\cal N}_Y\to{\cal N}_X$ is an isomorphism. For a monoid $N$ we denote by $N^{\rm gp}$ the associated group. For a finitely generated integral monoid $Q$ we let $$W[Q]=W\times_{\spec(\mathbb{Z})}\spec(\mathbb{Z}[Q])$$and give it the canonical log structure for which $Q$ is a chart.\\

\addtocounter{satz}{1}{\bf \arabic{section}.\arabic{satz} Definition} 
(i) A morphism of log schemes $f:(X,{\cal N}_X)\to(Y,{\cal N}_{{Y}})$ factors over the {\it log schematic image} $(\overline{X},{\cal N}_{\overline{X}})$ of $f$ which is defined as follows: The underlying scheme $\overline{X}$ is the schematic image of the morphism of schemes $X\to Y$ underlying $f$. Let $X\stackrel{i}{\to}\overline{X}\stackrel{\overline{f}}{\to} Y$ be the corresponding morphisms of schemes. The log structure ${\cal N}_{\overline{X}}$ is by definition the image of the natural composite map of log structures $\overline{f}^*{\cal N}_Y\to i_*f^*{\cal N}_Y\to i_*{\cal N}_X$ on $\overline{X}$. Here $i_*$ denotes the functor {\it push forward log structure}.\\(ii) A morphism of log schemes $f:(X,{\cal N}_X)\to(Y,{\cal N}_{{Y}})$ is said to be {\it log schematically dominant} if $({Y},{\cal N}_{{Y}})$ coincides with the log schematic image of $f$; it is said to be {\it log schematically dense} if in addition the underlying morphism of schemes is an open immersion.\\

A morphism of log schemes $i:(X,{\cal N}_X)\to(\overline{X},{\cal N}_{\overline{X}})$ is log schematically dense if and only if the underlying morphism of schemes is a schematically dense open immersion and the canonical morphism of log structures ${\cal N}_{\overline{X}}\to i_*{\cal N}_X$ is injective.\\ 

\begin{lem}\label{pfeil} Let $(X,{\cal N}_X)$ be a log scheme and $i:X\to\overline{X}$ a schematically dense open immersion of its underlying scheme into another scheme $\overline{X}$. Denote by $i_{*,{\rm sh}}{\cal N}_X$ the {\it sheaf theoretic push forward} of the sheaf of monoids ${\cal N}_X$. There exists a unique map $i_{*,{\rm sh}}{\cal N}_X\to (i_*{\cal N}_X)^{\rm gp}$ compatible with the natural maps $i_*{\cal N}_X\to i_{*,{\rm sh}}{\cal N}_X$ and $i_*{\cal N}_X\to(i_*{\cal N}_X)^{\rm gp}$.
\end{lem}

{\sc Proof:} First observe that ${\cal O}_{\overline{X}}\to i_*{\cal O}_X$ is injective, so henceforth we regard ${\cal O}_{\overline{X}}$ as a subsheaf of $i_*{\cal O}_X$. Also note $(i_*{\cal O}_X)^{\times}=i_*({\cal O}_X^{\times})$. It follows that we can view $i_*{\cal N}_X$ as the subsheaf of $i_{*,{\rm sh}}{\cal N}_X$ formed by those sections which map to ${\cal O}_{\overline{X}}$ under the map $\alpha:i_{*,{\rm sh}}{\cal N}_X\to i_*{\cal O}_X$ which we get by functoriality of sheaf theoretic push forward. To prove the lemma it is enough to show that $i_{*,{\rm sh}}{\cal N}_X$ arises from $i_*{\cal N}_X$ by inverting those sections $m$ for which the restrictions $\alpha(m)|_X$ are invertible. But this is the case: Take $m\in i_{*,{\rm sh}}{\cal N}_X$. Since $i_*{\cal O}_X$ arises from ${\cal O}_{\overline{X}}$ by inverting those Sections for which the restrictions to $X$ are invertible, we find $f,g\in{\cal O}_{\overline{X}}$ with $g|_X$ invertible and with $\alpha(m)=g^{-1}f$. We saw $g=\alpha(n)$ for some $n\in i_*{\cal N}_X$. Now $nm\in i_*{\cal N}_X$ and our claim and hence the lemma follows.\\

\begin{lem}\label{feinbild} The log schematic image $(\overline{X},{\cal N}_{\overline{X}})$ of a morphism of fine log schemes $f:(X,{\cal N}_X)\to(Y,{\cal N}_{{Y}})$ is a fine log scheme.
\end{lem}

{\sc Proof:} The coherence of ${\cal N}_{\overline{X}}$ follows from that of ${\cal N}_{{Y}}$. We have ${\cal N}_{\overline{X}}\subset i_*{\cal N}_X\subset i_{*,{\rm sh}}{\cal N}_X$, for the second inclusion see the proof of Lemma \ref{pfeil}. Therefore the integrality of ${\cal N}_X$ implies that of ${\cal N}_{\overline{X}}$.\\ 

\addtocounter{satz}{1}{\bf \arabic{section}.\arabic{satz} Definition}\newcounter{defbs1}\newcounter{defbs2}\setcounter{defbs1}{\value{section}}\setcounter{defbs2}{\value{satz}} A {\it log scheme with boundary} is a triple $((X,{\cal N}_X), (\overline{X},{\cal N}_{\overline{X}}),i)$ where $i:(X,{\cal N}_X)\to(\overline{X},{\cal N}_{\overline{X}})$ is a log schematically dense morphism such that $i^*{\cal N}_{\overline{X}}={\cal N}_X$ and $(i_*{\cal N}_X)^{\rm gp}={\cal N}_{\overline{X}}^{\rm gp}$. Let $(T,{\cal N}_T)$ be a log scheme. A $(T,{\cal N}_T)${\it -log scheme with boundary} is a log scheme with boundary $((X,{\cal N}_X), (\overline{X},{\cal N}_{\overline{X}}),i)$ together with a morphism of log schemes $g:(X,{\cal N}_X)\to(T,{\cal N}_T)$.\\ 

We think of $\overline{X}-X$ as a boundary of $X$. We will often drop $i$, $g$ and the log structures from our notation and just speak of the $T$-log scheme with boundary $(\overline{X},X)$. So in the following definition which justifies the whole concept.\\ 

\addtocounter{satz}{1}{\bf \arabic{section}.\arabic{satz} Definition:} The sheaf of relative differentials of a $T$-log scheme with boundary $(\overline{X},X)$ is defined as follows: Denote by $\tau$ the composite map$$i_{*,{\rm sh}}g^{-1}{\cal N}_T\to i_{*,{\rm sh}}{\cal N}_X\to (i_*{\cal N}_X)^{\rm gp}={\cal N}_{\overline{X}}^{\rm gp}$$where the second arrow is the one from Lemma \ref{pfeil}. Let $\Omega^1_{\overline{X}/W}$ be the sheaf of differentials of the morphism of underlying schemes $\overline{X}\to W$. Then $\Omega^1_{(\overline{X},X)/T}$ is the quotient of $$\Omega^1_{\overline{X}/W}\oplus({\cal O}_{\overline{X}}\otimes_{\mathbb{Z}}{\cal N}^{\rm gp}_{\overline{X}})$$ divided by the ${\cal O}_{\overline{X}}$-submodule generated by local sections of the forms
\begin{align}(d\alpha(a),0)-(0,\alpha(a)\otimes a)&\quad\quad\mbox{ with }a\in{\cal N}_{\overline{X}}\notag\\(0,1\otimes a)&\quad\quad\mbox{ with 
}a\in \bi(\tau).\notag \end{align}We define the de Rham complex $\Omega^{\bullet}_{(\overline{X},X)/T}$ by taking exterior powers and the differential as usual.\\

\begin{lem}\label{clabo} Let $(\overline{X},X)$ be a $T$-log scheme with boundary.\\(1) The restriction $\Omega^1_{(\overline{X},X)/T}|_X$ naturally coincides with the usual sheaf of relative logarithmic differentials of $g:(X,{\cal N}_X)\to(T,{\cal N}_T)$.\\(2) Suppose $g$ extends to a morphism of log schemes $\overline{g}:(\overline{X},{\cal N}_{\overline{X}})\to(T,{\cal N}_T)$. Let us assume the following conditions:
\begin{description}\item[(i)] The underlying scheme of $T$ is the spectrum of a field.\item[(ii)] For any \'{e}tale morphism $\overline{V}\to\overline{X}$ with $\overline{V}$ connected, the scheme $V=\overline{V}\times_{\overline{X}}X$ is also connected.\end{description}Then $\Omega^1_{(\overline{X},X)/T}$ naturally coincides with the usual sheaf $\Omega^1_{\overline{X}/T}$ of relative logarithmic differentials of $\overline{g}$.
\end{lem}

{\sc Proof:} (1) is immediate. (2) and its proof were suggested by the referee. Write $T=\spec(k)$. By base change, we may assume that $k$ is separably closed. It suffices to prove that the morphism $\overline{g}^{-1}{\mathcal N}_T\to i_{*,{\rm sh}}g^{-1}{\mathcal N}_T$ is an isomorphism. Let $x$ be a geometric point of $\overline{X}$ and let $\overline{V}$ be the strict Henselization of $\overline{X}$ at $x$. Put $V=\overline{V}\times_{\overline{X}}X$. Then, by the assumption (i), we see that both $\overline{V}$ and $V$ are connected. Hence we have$$(\overline{g}^{-1}{\mathcal N}_T)_x=\Gamma(\overline{V},\overline{g}^{-1}{\mathcal N}_T)=\Gamma(T,{\mathcal N}_T)$$$$(i_{*,{\rm sh}}g^{-1}{\mathcal N}_T)_x=\Gamma(\overline{V},\overline{g}^{-1}{\mathcal N}_T)=\Gamma(T,{\mathcal N}_T)$$and the lemma follows.\\

One class of examples where the condition (i) + (ii) of Lemma \ref{clabo} (2) holds true are the semistable $T$-log schemes discussed in Section 3; but for them, the conclusion $\Omega^1_{(\overline{X},X)/T}=\Omega^1_{\overline{X}/T}$ (if $g$ extends to $\overline{g}$) is immediate anyway. Undoubtly, if $g$ extends to $\overline{g}$, the conclusion of Lemma \ref{clabo} (2) holds under much more general conditions than the stated condition (i) + (ii). \\

\addtocounter{satz}{1}{\bf \arabic{section}.\arabic{satz} Examples:}\newcounter{bspcol1}\newcounter{bspcol2}\setcounter{bspcol1}{\value{section}}\setcounter{bspcol2}{\value{satz}} The following examples will be discussed later on.\\(a) Let $Q,P$ be finitely generated monoids and let $\rho:Q\to P^{\rm gp}$ be a morphism. Let $P'$ be the submonoid of $P^{\rm gp}$ generated by $P$ and $\bi(\rho)$. Then$$(W[P],W[P'])$$is a $T=W[Q]$-log scheme with boundary.\\(b1) Let $Q=\mathbb{N}$ with generator $t\in Q$. Let $t_1,\ldots,t_r$ be the standard generators of $\mathbb{N}^r$. Let $X=W[\mathbb{N}^r]$, the affine $r$-space over $W$ with the log structure defined by the divisor $V(t_1\cdot\ldots\cdot t_r)$. By means of $t\mapsto t_1\cdot\ldots\cdot t_r$ this is a $T=W[Q]$-log scheme. We compactify ${X}$ by$$\overline{X}=W\times_{\spec(\mathbb{Z})}(\times_{\spec(\mathbb{Z})}(\proj(\mathbb{Z}[t_0,t_i])_{1\le i\le r}))=({\bf P}_W^1)^r$$and take for ${\cal N}_{\overline{X}}$ the log structure defined by the normal crossing divisor $$(\overline{X}-X)\cup(\mbox{the closure of }V(t_1\cdot\ldots\cdot t_r)\subset X\mbox{ in }\overline{X}).$$(b2) Let $X$ and $T$ be as in (b1). Another compactifiction of $X$ is projective $r$-space, i.e. $\overline{X}'={\bf P}_W^r$; similarly we take ${\cal N}_{\overline{X'}}$ as the log structure defined by the normal crossing divisor $(\overline{X}'-X)\cup(\mbox{the closure of }V(t_1\cdot\ldots\cdot t_r)\subset X\mbox{ in }\overline{X}')$.\\(c) Let $k$ be a field, $W=\spec(k)$ and let again $Q=\mathbb{N}$ with generator $t\in Q$. The following type of $S=W[Q]$-log scheme with boundary (which generalizes \arabic{bspcol1}.\arabic{bspcol2}(b1) if $W=\spec(k)$ there) gives rise, by base change $t\mapsto 0$, to the $T$-log schemes with boundary discussed in Section \ref{semist} below. Let $\overline{X}$ be a smooth $W$-scheme, $X\subset \overline{X}$ a dense open subscheme, $D=\overline{X}-X$. Let $X\to S$ be a flat morphism, smooth away from the origin. Let $X_0$ be the fibre above the origin, let $\overline{X}_0$ be its schematic closure in $\overline{X}$ and suppose that $D\cup\overline{X}_0$ is a divisor with normal crossings on $\overline{X}$.\\
(d) Let $k$ be a field and let $T=(\spec(k),\mathbb{N}\stackrel{0}{\to}k)$, the standard logarithmic point (\cite{fkato}). Let $Y$ be a semistable $k$-log scheme in the sense of \cite{mokr} 2.4.1 or \cite{fkato}. That is, $Y$ is a fine $T$-log scheme $(Y,{\cal N}_Y)$ satisfying the following conditions. \'{E}tale locally on $Y$ there exist integers $i\ge 1$ and charts $\mathbb{N}^i\to{\cal N}_Y(Y)$ for ${\cal N}_Y$ such that\\(i) if on the log scheme $T$ we use the chart $\mathbb{N}\to k, 1\mapsto 0$, the diagonal morphism $\mathbb{N}\stackrel{\delta}{\to}\mathbb{N}^i$ is a chart for the structure morphism of log schemes $Y\to T$, and\\(ii) the induced morphism of schemes $$Y\to\spec(k)\times_{\spec(k[t])}\spec(k[t_1,\ldots,t_i])$$ is smooth in the classical sense. Let $\overline{X}$ be the union of some irreducible components of $Y$ and let $X$ be the open subscheme of $\overline{X}$ which is the complement in $Y$ of the union of all irreducible components not contained in $\overline{X}$. Then $\overline{X}$ inherits a structure of $T$-log scheme, but it is not log smooth over $T$. However, we can view $(\overline{X},X)$ as a $T$-log scheme with boundary (forgetting that the morphism $X\to T$ actually extends to $\overline{X}$): as such it is what we will call {\it smooth} below.\\

\addtocounter{satz}{1}{\bf \arabic{section}.\arabic{satz}}\newcounter{bspbdl1}\newcounter{bspbdl2}\setcounter{bspbdl1}{\value{section}}\setcounter{bspbdl2}{\value{satz}} A concrete example for \arabic{bspcol1}.\arabic{bspcol2}(c) (see \cite{hkstrat} for more details). Again let $k$ be a field and let $S=W[\mathbb{N}]$ with generator $q$ of $\mathbb{N}$. Let $Y$ be a semistable $k$-log scheme with set of irreducible components $\{Y_j\}_{j\in R}$ all of which we assume to be smooth. As in \cite{kalo} p.222/223 we define for every $j\in R$ an invertible ${\cal O}_Y$-module ${\cal F}_j$ as follows: Let ${\cal N}_{Y,j}$ be the subsheaf of the log structure ${\cal N}_Y$ of $Y$ which is the preimage of $\ke({\cal O}_Y\to{\cal O}_{Y_j})$. This ${\cal N}_{Y,j}$ is a principal homogeneous space over ${\cal O}_Y^{\times}$, and its associated invertible ${\cal O}_Y$-module is ${\cal F}_j$. Now fix a subset $I\subset R$ with $|I|=i$ and let $L=R-I$. Suppose $M=\cap_{j\in I}Y_j$ is nonempty. Let $$V_M=\spec(\sym_{{\cal O}_M}(\oplus({\cal F}_j)_{j\in I}))=\times_M(\spec(\sym_{{\cal O}_M}({\cal F}_j)))_{j\in I}.$$By its definition, the affine vector bundle $V_M$ over $M$ comes with a natural coordinate cross, a normal crossing divisor on $V_M$. The intersection of $M$ with all irreducible components of $Y$ not containing $M$ is a normal crossing divisor $D$ on $M$. Let $D_V'\subset V_M$ be its preimage under the structure map $V_M\to M$ and let $D_V\subset V_M$ be the union of $D_V'$ with the natural coordinate cross in $V_M$. Then $D_V$ is a normal crossing divisor on $V_M$. Let ${\cal N}_{V_M}$ be the corresponding log structure on $V_M$. There exists a distinguished element $a\in\Gamma(V_M,{\cal O}_M)$ having $D_V$ as its set of zeros and such that the assignment $q\mapsto a$ defines a morphism of log schemes $V_M\to S$ with the following property: The induced $S$-log scheme $(M,{\cal N}_{V_M}|_M)$ on the zero Section $M\to V_M$ coincides with the $S$-log scheme $(M,{\cal N}_Y|M)$ induced by $Y$. This $a\in\Gamma(V_M,{\cal O}_M)=\sym_{{\cal O}_M}(\oplus({\cal F}_j)_{j\in I})(M)$ can be described as follows: Denote the image of $q\in{\cal N}_S(S)$ (here ${\cal N}_S$ is the log structure of $S$) under the structure map ${\cal N}_S(S)\to {\cal N}_Y(Y)\to {\cal N}_Y|_M(M)$ again by $q$. Locally on $M$ it can be (non-uniquely) factored as $q=t_0\prod_{j\in I}v_j$ where $v_j$ is a (local) generator of ${\cal F}_j|_M$ and $t_0$ maps to a (local) defining equation $a_0\in{\cal O}_M$ of the divisor $D$ in $M$. Then $a=a_0.(\oplus_{j\in I}v_j)\in\sym_{{\cal O}_M}(\oplus_{j\in I}{\cal F}_j)(M)$ is the wanted element, globally well defined. We can view $V_M$ in a canonical way as a (schematically) dense open subscheme of $$P_M=\times_M(\proj(\sym_{{\cal O}_M}({\cal O}_M\oplus{\cal F}_j)))_{j\in I}$$by identifying a homogenous section  $s\in\sym_{{\cal O}_M}({\cal F}_j)$ of degree $n$ with the degree zero section $s/1_{{\cal O}_M}^n$ of $\sym_{{\cal O}_M}({\cal O}_M\oplus{\cal F}_j)[1_{{\cal O}_M}^{-1}]$. We give $P_M$ the log structure defined by the normal crossing divisor $(P_M-V_M)\cup \overline{D}_V$, where $\overline{D}_V$ is the closure of $D_V$ in $P_M$. Then $(P_M,V_M)$ is a $S$-log scheme with boundary.\\ 

\addtocounter{satz}{1}{\bf \arabic{section}.\arabic{satz}}\newcounter{basech1}\newcounter{basech2}\setcounter{basech1}{\value{section}}\setcounter{basech2}{\value{satz}} A {\it morphism} of $T$-log schemes with boundary $f:(\overline{X},X)\to(\overline{X}',X')$ is a morphism of log schemes $$f:(\overline{X},{\cal N}_{\overline{X}})\to(\overline{X}',{\cal N}_{\overline{X}'})$$ with $X\subset f^{-1}(X')$ and restricting to a morphism of $T$-log schemes $(X,{\cal N}_{{X}})\to({X}',{\cal N}_{{X}'})$. We have a fully faithful functor from the category of $T$-log schemes to the category of $T$-log schemes with boundary. Namely, take $Y$ to $(Y,Y)$. Beware that $(T,T)$ is {\it not} a final object in the category of $T$-log schemes with boundary. We have obvious base change functors for morphisms $W'\to W$ to our underlying base scheme $W$ and everything we develop here behaves well with respect to these base changes. We also have {\it base change functors for closed immersions of log schemes} $T'\to T$ as follows: if $(\overline{X},X)$ is a $T$-log scheme with boundary, let $X_{T'}=X\times_TT'$ be the fibre product in the category of log schemes. Define the log scheme $\overline{X}_{T'}$ as the log schematic image of the morphism of log schemes $X_{T'}\to\overline{X}$. Then $(\overline{X}_{T'},X_{T'})$ is a $T'$-log scheme with boundary.\\ 

\addtocounter{satz}{1}{\bf \arabic{section}.\arabic{satz}} For the rest of this paper we always assume that the log scheme $T$ is fine. All fibre products of fine log schemes are taken in the category of fine log schemes, unless specified otherwise. A $T$-log scheme with boundary $(\overline{X},X)$ is said to be {\it fine} if the log scheme $(\overline{X},{\cal N}_{\overline{X}})$ is fine.

\begin{lem} In the category of fine $T$-log schemes with boundary, products  exist.
\end{lem}

{\sc Proof:} Given fine $T$-log schemes with boundary $(\overline{X}_1,X_1)$ and $(\overline{X}_2,X_2)$, set $$(\overline{X}_1,X_1)\times_T(\overline{X}_2,X_2)=(\overline{X}_1\overline{\times}_T\overline{X}_2,X_1\times_TX_2).$$Here $X_1\times_TX_2$ denotes the fibre product in the category of fine $T$-log schemes, and $\overline{X}_1\overline{\times}_T\overline{X}_2$ is defined as the log schematic image of $X_1\times_TX_2\to \overline{X}_1\times_{W}\overline{X}_2$. (So $\overline{X}_1\overline{\times}_T\overline{X}_2$ depends also on $X_1$ and $X_2$, contrary to what the notation suggests. Note that by the construction \cite{kalo} 2.7, the scheme underlying $X_1\times_TX_2$ is a subscheme of the scheme theoretic fibre product, hence is a subscheme of the scheme underlying $\overline{X}_1\times_{W}\overline{X}_2$.) That $\overline{X}_1\overline{\times}_T\overline{X}_2$ is fine follows from Lemma \ref{feinbild}.\\

\addtocounter{satz}{1}{\bf \arabic{section}.\arabic{satz}}\newcounter{gegenb1}\newcounter{gegenb2}\setcounter{gegenb1}{\value{section}}\setcounter{gegenb2}{\value{satz}} It is to have fibre products why we did not require $X=f^{-1}(X')$ in the definition of morphisms of $T$-log schemes with boundary. If the structural map from the underlying scheme of the log scheme $T$ to $W$ is an isomorphism, one has $(\overline{X},X)\cong(\overline{X},X)\times_T(T,T)$. However, we stress that in contrast to taking the base change with the identity $T\to T$ (cf. \arabic{basech1}.\arabic{basech2}), the operation of taking the fibre product with the $T$-log scheme with boundary $(T,T)$ is non-trivial in general. For example, let $Q=\mathbb{N}$ with generator $q\in Q$, let $T=W[Q]$ and let $U_1, U_2$ be the standard generators of $\mathbb{N}^2$. For $i\in\mathbb{Z}$ let $\overline{X}_i=W[\mathbb{N}^2]$, and let $X_i=W[\mathbb{Z}\oplus\mathbb{N}]$, the open subscheme of $\overline{X}_i$ where $U_1$ is invertible. Define a structure of $T$-log scheme with boundary on $(\overline{X}_i,X_i)$ by sending $q\mapsto U_1^iU_2$. Then $$(\overline{X}_i,X_i)\cong(\overline{X}_i,X_i)\times_T(T,T)\quad\mbox{ if }i\ge0$$$$(\overline{X}_i,X_i)\not\cong(\overline{X}_i,X_i)\times_T(T,T)\quad\mbox{ if }i<0.$$Indeed, $\overline{X}_i\overline{\times}_TT$ is the closure in $W[Q\oplus\mathbb{N}^2]$ of the closed subscheme $V(q-U_1^iU_2)$ of $W[Q\oplus\mathbb{Z}\oplus\mathbb{N}]$. If $i\ge0$ this is the subscheme $V(q-U_1^iU_2)$ of $W[Q\oplus\mathbb{N}^2]$ which maps isomorphically to $W[\mathbb{N}^2]$. If $i<0$ this is the subscheme $V(qU_1^{-i}-U_2)$ of $W[Q\oplus\mathbb{N}^2]$ which does not map isomorphically to $W[\mathbb{N}^2]$.\\

\addtocounter{satz}{1}{\bf \arabic{section}.\arabic{satz}}\newcounter{defcha1}\newcounter{defcha2}\setcounter{defcha1}{\value{section}}\setcounter{defcha2}{\value{satz}} Let $(\overline{X},X)$ be a fine $T$-log scheme with boundary. A {\it chart} $(Q\to P^{\rm gp}\supset P)$ for $(\overline{X},X)$ over $T$ is a chart $\lambda:P\to\Gamma(\overline{X},{\cal N}_{\overline{X}})$ for $(\overline{X},{\cal N}_{\overline{X}})$, a chart $\sigma:Q\to \Gamma({T},{\cal N}_{{T}})$ for $({T},{\cal N}_{{T}})$ and a morphism $\rho:Q\to P^{\rm gp}$ such that $\lambda^{\rm gp}\circ\rho=\tau\circ\sigma$, where $\tau:\Gamma({T},{\cal N}_{{T}})\to\Gamma({X},{\cal N}_{{X}})\to\Gamma(\overline{X},{\cal N}_{\overline{X}}^{\rm gp})$ is the composite of the structural map with that from Lemma \ref{pfeil}. 

\begin{lem} \'{E}tale locally on $\overline{X}$, charts for $(\overline{X},X)$ exist.
\end{lem}

{\sc Proof:} (corrected version due to the referee) We may by \cite{kalo} assume that $(\overline{X},{\cal N}_{\overline{X}})$ has a chart $g:G\to \Gamma(\overline{X},{\cal N}_{\overline{X}})$ and $({T},{\cal N}_{{T}})$ has a chart $\sigma:Q\to \Gamma({T},{\cal N}_{{T}})$. Let $x\in X$ and let ${\cal N}_{\overline{X},\overline{x}}$ be the stalk of ${\cal N}_{\overline{X}}$ at the separable closure $\overline{x}$ of $x$. Let $\varphi$ be the composite $$Q\stackrel{\sigma}{\to}\Gamma({T},{\cal N}_{{T}})\stackrel{\tau}{\to}\Gamma(\overline{X},{\cal N}_{\overline{X}}^{\rm gp})\to{\cal N}_{\overline{X},\overline{x}}^{\rm gp}.$$ Choose generators $q_1,\ldots,q_m$ of $Q$ and elements $x_i, y_i\in {\cal N}_{\overline{X},\overline{x}}$ $(1\le i\le m)$ such that $\varphi(q_i)=x_iy_i^{-1}$. Next, choose elements $a_i, b_i\in G$ and $u_i, v_i\in {\cal O}^{\times}_{\overline{X},\overline{x}}$ $(1\le i\le m)$ satisfying $g(a_i)=x_iu_i$ and $g(b_i)=y_iv_i$: these elements exist because $g$ is a chart. Now let$$f:G^{\rm gp}\oplus Q^{\rm gp}\oplus\mathbb{Z}^m\oplus\mathbb{Z}^m\longrightarrow{\cal N}_{\overline{X},\overline{x}}^{\rm gp}$$be the morphism defined by$$(h,q,(k_i)_{i=1}^m,(l_i)_{i=1}^m)\mapsto g^{\rm gp}(h)\varphi^{\rm gp}(q)\prod_{i=1}^mu_i^{k_i}\prod_{i=1}^mv_i^{l_i},$$and define $P$ by $P=f^{-1}({\cal N}_{\overline{X},\overline{x}})$. Then $f|_P:P\to {\cal N}_{\overline{X},\overline{x}}$ extends to a chart around $\bar{x}$ by \cite{kalo} 2.10. It remains to prove that the canonical inclusion $Q\to G^{\rm gp}\oplus Q^{\rm gp}\oplus\mathbb{Z}^m\oplus\mathbb{Z}^m, q\mapsto (1,q,0,0)$ actually takes values in $P^{\rm gp}$. Write a given $q\in Q$ as $q=\prod_{i=1}^mq_i^{n_i}$ with $n_i\in\mathbb{N}$. Then we have

$$f(q)p=\prod_{i=1}^m(\frac{x_i}{y_i})^{n_i}=\prod_{i=1}^m(\frac{x_iu_i}{y_iv_i}\cdot\frac{v_i}{u_i})^{n_i}=\frac{f((\prod_ia_i^{n_i},0,(0),(n_i)_i))}{f((\prod_ib_i^{n_i},0,(n_i)_i,(0)))}.$$
Put $\alpha=(\prod_ia_i^{n_i},0,(0),(n_i)_i)$ and $\beta=(\prod_ib_i^{n_i},0,(n_i)_i,(0))$. Then we have $\alpha, \beta\in P$ and $f(q\beta)=f(\alpha)$. So $q\beta$ is in $P$ by the definition of $P$ and so $q$ maps to $P^{\rm gp}$.\\

\section{Smoothness}

\addtocounter{satz}{1}{\bf \arabic{section}.\arabic{satz} Definition:}\newcounter{defsmo1}\newcounter{defsmo2}\setcounter{defsmo1}{\value{section}}\setcounter{defsmo2}{\value{satz}} (1) A morphism of $T$-log schemes with boundary $(\overline{Y},Y)\to(\overline{X},X)$ is said to be a {\it boundary exact closed immersion} if $\overline{Y}\to\overline{X}$ is an exact closed immersion and if for every open neighbourhood $U$ of $Y$ in $X$, there exists an open neighbourhood $\overline{U}$ of $\overline{Y}$ in $\overline{X}$ with $U$ schematically dense in $\overline{U}$.\\ (2) A {\it first order thickening} of $T$-log schemes with boundary is a morphism $(\overline{L}',L')\to(\overline{L},L)$ such that $\overline{L'}\to\overline{L}$ is an exact closed immersion defined by a square zero ideal in ${\cal O}_{\overline{L}}$.\\(3) A fine $T$-log scheme with boundary $(\overline{X},X)$ is said to be {\it weakly smooth} if $\overline{X}$ is locally of finite presentation over $W$ and if the following condition holds: for every first order thickening $\eta:(\overline{L}',L')\to(\overline{L},L)$ and for every morphism $\mu:(\overline{L}',L')\to(\overline{X},X)$ there is \'{e}tale locally on $\overline{L}$ a morphism $\epsilon:(\overline{L},L)\to(\overline{X},X)$ such that $\mu=\epsilon\circ\eta$.\\(4) A $T$-log scheme with boundary $(\overline{X},X)$ is said to be {\it smooth} if it is weakly smooth and satisfies the following property: For all morphisms $(\overline{Y},Y)\to(\overline{X},X)$ and all boundary exact closed immersions $(\overline{Y},Y)\to(\overline{V},V)$ of fine $T$-log schemes with boundary, there exists \'{e}tale locally on $(\overline{X}\overline{\times}_T\overline{V})$ an exactification $$\overline{Y}\to Z\to(\overline{X}\overline{\times}_T\overline{V})$$of the diagonal embedding $\overline{Y}\to(\overline{X}\overline{\times}_T\overline{V})$ (a morphism of log schemes in the usual sense) such that the projection $Z\to(\overline{X}\overline{\times}_T\overline{V})\to \overline{V}$ is strict and log smooth.\\

Recall that by \cite{kalo} 3.8, `strict and log smooth' is equivalent to `strict and smooth on underlying schemes'. A $T$-log scheme $X$ is log smooth if and only if $(X,X)/T$ is weakly smooth. Assume this is the case. Then $(X,X)/T$ satisfies the smoothness condition with respect to test objects $({X},X){\leftarrow}(\overline{Y},Y)\rightarrow({V},V)$ (i.e. for which $\overline{V}=V$), because ${X}\overline{\times}_T{V}\stackrel{p}{\to}{V}$ is clearly log smooth. For general $(\overline{V},V)$ (and log smooth $T$-log schemes $X$) we have at least Theorem \ref{sglatt} and Theorem \ref{tglatt} below (note that the hypotheses of Proposition \ref{wecri} below for $(X,X)/T$ are {\it equivalent} to log smoothness of $X/T$, by \cite{kalo} 3.5 and as worked out in \cite{fkato}).\\ 

\begin{pro}\label{smocri} Let $(\overline{X},X)$ be a weakly smooth $T$-log scheme with boundary and let $T_1\to T$ be an exact closed immersion. Then $(\overline{X}_{T_1},X_{T_1})$ is a weakly smooth $T_1$-log scheme with boundary.
\end{pro}

{\sc Proof:} Let $$(\overline{X}_{T_1},X_{T_1})\stackrel{\mu}{\leftarrow}(\overline{L}',L')\stackrel{\eta}{\to}(\overline{L},L)$$ be a test object over $T_1$. By the weak smoothness of $(\overline{X},X)/T$ we get $\epsilon:(\overline{L},L)\to(\overline{X},X)$ \'{e}tale locally on $\overline{L}$ such that $\mu=\epsilon\circ\eta$. The restriction $\epsilon|_L:L\to\overline{X}$ goes through $X_{T_1}$; since $L$ is log schematically dense in $\overline{L}$ this implies that $\epsilon$ goes through $(\overline{X}_{T_1},X_{T_1})$ (the schematic image is transitive, \cite{EGA} I, 9.5.5).\\

\begin{pro}\label{wecri} Suppose $W$ is locally noetherian. Let $Q$ be a finitely generated integral monoid, let $S=W[Q]$ and let $T\to S$ be an exact closed immersion. Let $(\overline{X},X)$ be a $T$-log scheme with boundary. Suppose that \'{e}tale locally on $\overline{X}$ there are charts $Q\to P^{\rm gp}\supset P$ for $(\overline{X},X)$ over $T$ as in \arabic{defcha1}.\arabic{defcha2} such that the following conditions ${\rm (i), (ii)}$ are satisfied:\begin{description}\item[(i)] The kernel and the torsion part of the cokernel of $Q^{\rm gp}\to P^{\rm gp}$ are finite groups of orders invertible on $W$.\item[(ii)] Let $P'$ be the submonoid of $P^{\rm gp}$ generated by $P$ and the image of $Q\to P^{\rm gp}$ and let $W[P]_T$ be the schematic closure of $W[P']\times_ST=W[P']_T$ in $W[P]$. Then $\lambda:\overline{X}\to W[P]_T$ is smooth on underlying schemes.\end{description}Then $(\overline{X},X)/T$ is weakly smooth.
\end{pro}

{\sc Proof:} (Note that $\lambda$ in (ii) exists by the schematic density of $X\to\overline{X}$.) Let $$(\overline{X},X)\stackrel{\mu}{\leftarrow}(\overline{L}',L')\stackrel{\eta}{\to}(\overline{L},L)$$ be a test object over $T$. Using (i), one can follow the arguments in \cite{kalo} 3.4 to construct morphisms $(\overline{L},L)\to(W[P],W[P'])$ of $S$-log schemes with boundary. Necessarily $L$ maps in fact to $W[P']_T$. Since $L\to \overline{L}$ is log schematically dense, $\overline{L}$ maps in fact to $W[P]_T$. By (ii) this morphism can be lifted further to a morphism $\overline{L}\to\overline{X}$ inducing $(\overline{L},L)\to(\overline{X},X)$ as desired.\\

\begin{satz}\label{sglatt} In the situation of Proposition \ref{wecri}, suppose in addition $S=T$ and $T\to S$ is the identity. Then for every $S$-log scheme with boundary $(\overline{V},V)$, the projection $\overline{X}\overline{\times}_S\overline{V}\stackrel{p}{\to}\overline{V}$ is log smooth.
\end{satz} 

{\sc Proof:} We may assume that $(\overline{X},X)$ over $T$ has a chart as described in Proposition \ref{wecri} and that $(\overline{V},V)$ over $T$ has a chart $Q\to F^{\rm gp}$, $F\to{\cal N}_{\overline{V}}(\overline{V})$. Our assumptions imply that$$\overline{X}\times_W\overline{V}\to W[P]\times_W\overline{V}$$is smooth on underlying schemes. It is also strict, hence log smooth. Perform the base change with the closed immersion of log schemes$$W[P]\overline{\times}_S\overline{V}\to W[P]\times_W\overline{V}$$to get the log smooth morphism$$\overline{X}\overline{\times}_S\overline{V}\stackrel{}{\to}W[P]\overline{\times}_S\overline{V}$$(by our construction of fibre products, $W[P]\overline{\times}_S\overline{V}$ is the log schematic closure of $W[P']{\times}_S{V}$). Its composite with the projection $$W[P]\overline{\times}_S\overline{V}\stackrel{\beta}{\to}\overline{V}$$is $p$, hence it is enough to show that $\beta$ is log smooth. Now $\beta$ arises by the base change $\overline{V}\to W[F]$ from the projection $$W[P]\overline{\times}_SW[F]\stackrel{\gamma}{\to}W[F]$$so that it is enough to show that $\gamma$ is log smooth. Let $F'$ be the submonoid of $F^{\rm gp}$ generated by $F$ and the image of $Q\to F^{\rm gp}$. Let $(P'\oplus_QF')^{\rm int}$ be the push out of $P'\leftarrow Q\to F'$ in the category of integral monoids, i.e. $(P'\oplus_QF')^{\rm int}=\bi(P'\oplus_QF'\to (P'\oplus_QF')^{\rm gp})$ where $P'\oplus_QF'$ is the push out in the category of monoids. (If $Q$ is generated by a single element then actually $(P'\oplus_QF')^{\rm int}=P'\oplus_QF'$ by \cite{kalo} 4.1.) Define the finitely generated integral monoid $$R=\bi(P\oplus F\to (P'\oplus_QF')^{\rm int}).$$Then $\gamma$ can be identified with the natural map $W[R]\to W[F]$. That this is log smooth follows from \cite{kalo} 3.4 once we know that $$a:F^{\rm gp}\to R^{\rm gp}=(P^{\rm gp}\oplus_{Q^{\rm gp}}F^{\rm gp})$$has kernel and torsion part of the cokernel finitely generated of orders invertible on $W$. But this follows from the corresponding facts for $b:Q^{\rm gp}\to P^{\rm gp}$ because we have isomorphisms $\ke(b)\cong\ke(a)$ and $\koke(b)\cong\koke(a)$.\\

\begin{satz}\label{tglatt} In the situation of Proposition \ref{wecri},  $(\overline{X},X)/T$ is smooth.
\end{satz}

{\sc Proof:} It remains to verify the second condition in the definition of smoothness. Let $(\overline{Y},Y)\to(\overline{X},X)$ and $(\overline{Y},Y)\to(\overline{V},V)$ be test objects. We may assume that $\overline{Y}$ is connected. Remove all irreducible components of $\overline{V}$ not meeting $\bi(\overline{Y})$ so that we may assume that each open neighbourhood of $\bi(\overline{Y})$ in $\overline{V}$ is schematically dense. After \'{e}tale localization we may assume that $(\overline{X},X)$ has a chart $P\to\Gamma(\overline{X},{\cal N}_{\overline{X}}^{\rm gp})$ as in Proposition \ref{wecri}. Viewing our test objects as objects over $S$ we can form the fibre product of fine $S$-log schemes with boundary $(W[P]\overline{\times}_S\overline{V},W[P']\times_SV)$. \'{E}tale locally on $W[P]\overline{\times}_S\overline{V}$ we find an exactification $$\overline{Y}\stackrel{i}{\to} \tilde{Z}\stackrel{\tilde{g}}{\to} W[P]\overline{\times}_S\overline{V}$$of the diagonal embedding $\overline{Y}\to W[P]\overline{\times}_S\overline{V}$. We may assume that $\tilde{Z}$ is connected. After further \'{e}tale localization on $\tilde{Z}$ we may also assume that $\tilde{q}=\tilde{p}\circ \tilde{g}:\tilde{Z}\to\overline{V}$ is strict, where $\tilde{p}:(W[P]\overline{\times}_S\overline{V})\to \overline{V}$ is the projection: this follows from the fact that for $y\in \overline{Y}$ the stalks of the log structures ${\cal N}_{\tilde{Z}}$ and $\tilde{q}^*{\cal N}_{\overline{V}}$ at the separable closure of $i(y)$ coincide, because $\overline{Y}\stackrel{i}{\to} \tilde{Z}$ and $\overline{Y}\to \overline{V}$ are exact closed immersions. By Theorem \ref{sglatt}, $\tilde{p}$ is log smooth. Thus $\tilde{q}$ is also log smooth, hence is smooth on underlying schemes. Let $$\tilde{Z}^0=\tilde{Z}\times_{(W[P]\overline{\times}_S\overline{V})}(W[P']\times_SV),$$an open subscheme of $\tilde{Z}$ containing $\bi(Y)$. Consider the restriction $\tilde{q}^0:\tilde{Z}^0\to\overline{V}$ of $\tilde{q}$. Since it is smooth on underlying schemes, it maps schematically dominantly to an open neighbourhood of $\bi(Y)$ in $V$ (here a morphism of schemes ${\cal X}\to{\cal Y}$ is said to be schematically dominant if its schematic image coincides with ${\cal Y}$). It follows that $\tilde{q}^0$ maps schematically dominantly also to $\overline{V}$ because of our assumption on $\overline{V}$ and the fact that $(\overline{Y},Y)\to(\overline{V},V)$ is boundary exact. Thus $\tilde{q}$ is a classically smooth morphism from the connected scheme $\tilde{Z}$ to another scheme $\overline{V}$ such that its restriction to the open subscheme $\tilde{Z}^0$ maps schematically dominantly to $\overline{V}$. This implies that $\tilde{Z}^0$ is schematically dense in $\tilde{Z}$, because (schematically) dominant classically smooth morphisms from a connected scheme induce bijections between the respective sets of irreducible components. It follows that $\tilde{g}$ factors as $$\tilde{Z}\stackrel{g}{\to}(W[P]_T\overline{\times}_T\overline{V})\stackrel{k}{\to}W[P]\overline{\times}_S\overline{V}:$$first as a morphism of underlying schemes because its restriction to the open schematically dense subscheme $\tilde{Z}^0$ factors through $$W[P']_T\times_T{V}=W[P']\times_S{V};$$but then also as a morphism of log schemes, because $k$ is strict. The morphism $g$ is log \'{e}tale because the composite $\tilde{g}$ with the closed embedding $k$ is log \'{e}tale. Let $$Z=\tilde{Z}\times_{(W[P]_T\overline{\times}_T\overline{V})}(\overline{X}\overline{\times}_T\overline{V}).$$From the assumption (ii) in Proposition \ref{wecri} we deduce that $\overline{X}\overline{\times}_T\overline{V}\to W[P]_T\overline{\times}_T\overline{V}$ is log smooth and strict, hence $Z\to\tilde{Z}$ is log smooth and strict, hence smooth on underlying schemes. Together with the smoothness of $\tilde{q}$ it follows that $Z\to\overline{V}$ is smooth on underlying schemes. Furthermore $Z\to \overline{X}\overline{\times}_T\overline{V}$ is log \'{e}tale because $g$ is log \'{e}tale. Finally, $Y\to Z$ is an exact closed immersion because $Z\to\tilde{Z}$ is strict and $Y\to\tilde{Z}$ is an exact closed immersion. The theorem is proven.\\

The interest in smoothness as we defined it lies in the following proposition, which enables us to develop nice cohomology theories for $T$-log schemes with boundary.\\

\begin{pro}\label{keykoh} Let $(\overline{Y},Y)\to(\overline{X}_i,X_i)$ be boundary exact closed immersions into smooth $T$-log schemes with boundary ($i=1,2$). Then there exist \'{e}tale locally on $(\overline{X}_1\overline{\times}_T\overline{X}_2)$ factorizations$$(\overline{Y},Y)\stackrel{\iota}{\to}(\overline{Z},Z)\stackrel{}{\to}(\overline{X}_1\overline{\times}_T\overline{X}_2,X_1\times_TX_2)$$of the diagonal embedding such that $\iota$ is a boundary exact closed immmersion, the map $\overline{Z}\to\overline{X}_1\overline{\times}_T\overline{X}_2$ is log \'{e}tale, and the projections $p_i:\overline{Z}\to\overline{X}_i$ are strict and log smooth, hence smooth on underlying schemes.
\end{pro}

{\sc Proof:} By the definition of smoothness we find \'{e}tale locally exactifications ($i=1,2$)$$\overline{Y}\to \overline{Z}_i\to\overline{X}_1\overline{\times}_T\overline{X}_2$$such that the projections $\overline{Z}_i\to \overline{X}_i$ are strict and log smooth. Let $$\overline{Z}'=\overline{Z}_1\times_{(\overline{X}_1\overline{\times}_T\overline{X}_2)}\overline{Z}_2$$and let $\overline{Y}\to\overline{Z}\to\overline{Z}'$ be an exactification of $\overline{Y}\to\overline{Z}'$. After perhaps \'{e}tale localization on $\overline{Z}$ as in the proof of Theorem \ref{tglatt} we may assume that the projections $\overline{Z}\to\overline{Z}_i$ are strict. Hence the projections $p_i:\overline{Z}\to\overline{X}_i$ are strict and log smooth. This implies that $$Z=p_1^{-1}(X_1)\cap p_2^{-1}(X_2)$$is log schematically dense in $\overline{Z}$. Indeed, it suffices to prove the log schematic density of $Z$ in $p_1^{-1}(X_1)$ and of $p_1^{-1}(X_1)$ in $\overline{Z}$. Both assertions follow from the general fact that for a strict and log smooth (and in particular classically smooth) morphism of log schemes $h:L\to S$ and a log schematically dense open immersion $S'\to S$, also $h^{-1}(S')$ with its pull back log structure from $S'$ is log schematically dense in $L$: this is easy to prove since the question is local for the \'{e}tale topology and we therefore may assume that $h$ is a relative affine space. The classical smoothness of (say) $p_1$ and the boundary exactness of $(\overline{Y},Y)\to(\overline{X}_1,X_1)$ imply that $(\overline{Y},Y)\to(\overline{Z},Z)$ is boundary exact (for each connnected component $\overline{Z}'$ of $\overline{Z}$ the map $\pi_0(\overline{Z}')\to\pi_0(\overline{X}_1)$ between sets of irreducible components induced by $p_1$ is injective). We are done.\\

\addtocounter{satz}{1}{\bf \arabic{section}.\arabic{satz} Examples:} We make the exactification $Z\to\overline{X}\overline{\times}_T\overline{V}$ in Theorem \ref{tglatt} explicit in some examples, underlining the delicacy of the base change argument in Theorem \ref{tglatt}. In the following, for free variables $U_1,\ldots,V_1,\ldots$ we denote by $W[U_1,\ldots,V_1^{\pm},\ldots]$ the log scheme $$W[\mathbb{N}\oplus\ldots\oplus\mathbb{Z}\oplus\ldots]$$ with generators $U_1,\ldots$ for $\mathbb{N}\oplus\ldots$ and generators $V_1,\ldots$ for $\mathbb{Z}\oplus\ldots$. For $f\in\mathbb{Z}[U_1,\ldots,V_1^{\pm},\ldots]$ we denote by $W[U_1,\ldots,V_1^{\pm},\ldots]/f$ the exact closed subscheme defined by $f$.\\(a) Let $Q=\mathbb{N}$ with generator $q$. Let $X=W[U_1^{\pm},U_2]\subset\overline{X}=W[U_1,U_2]$. Define $X\to S$ by sending $q\mapsto U_1^{-1}U_2$, thus $(\overline{X},X)$ is a smooth $S$-log scheme with boundary. The self fibre product of $S$-log schemes with boundary is $$(\overline{X}_1,X_1)=(\overline{X},X)\overline{\times}_S(\overline{X},X)\quad\quad\quad\quad\quad\quad\quad\quad\quad\quad\quad\quad\quad\quad\quad\quad\quad\quad\quad\quad\quad\quad\quad\quad\quad$$$$=(W[U_1,U_2,V_1,V_2]/(V_1U_2-V_2U_1),W[U_1^{\pm},U_2,V_1^{\pm},V_2]/(U_1^{-1}U_2-V_1^{-1}V_2)).$$Note that the projections $q_j:\overline{X}_1\to\overline{X}$ are not flat (the fibres above the respective origins are two dimensional), although they are log smooth. We construct the desired log \'{e}tale map $Z\stackrel{g}{\to}\overline{X}_1$ according to the procedure in \cite{kalo}, 4.10. Embed $\mathbb{Z}\to\mathbb{Z}^4$ by sending $n\mapsto(n,-n,-n,n)$ and let $H$ be the image of the canonical map $\mathbb{N}^4\to(\mathbb{Z}^4/\mathbb{Z})$. Then $\overline{X}_1=W[H]$. Let $h:(\mathbb{Z}^4/\mathbb{Z})\to\mathbb{Z}^2$ be the map which sends the class of $(n_1,n_2,n_3,n_4)$ to $(n_1+n_3,n_2+n_4)$, and let $K=h^{-1}(\mathbb{N}^2)$. Then $Z=W[K]$ works. More explicitly: We have an isomorphism $K\cong\mathbb{N}^2\oplus\mathbb{Z}$ by sending the class of $(n_1,n_2,n_3,n_4)$ to $(n_1+n_3,n_2+n_4,n_1+n_2)$. Then $$Z=W[S_1,S_2,S_3^{\pm}]$$ and $g$ is given by $U_1\mapsto S_1S_3,\quad U_2\mapsto S_2S_3,\quad V_1\mapsto S_1,\quad V_2\mapsto S_2$.

Now consider the base change with $T=W[q]/q\to S$ defined by sending $q\mapsto 0$. For $j=1,2$ let $\overline{X}_{1,j}=\overline{X}_1\times_{\overline{X}}\overline{X}_T$ where in the fibre product we use the $j$-th projection as the structure map for the first factor. Let $\overline{X}_{T,1}=\overline{X}_T\overline{\times}_T\overline{X}_T$. Then we find $\overline{X}_{1,1}=W[U_1,V_1,V_2]/(V_2U_1)$, $\overline{X}_{1,2}=W[U_1,U_2,V_1]/(V_1U_2)$, thus containing $\overline{X}_{T,1}=W[U_1,V_1]$ as a {\it proper} subscheme.\\(b) Let $S, X, \overline{X}$ be as in (a), but this time define $X\to S$ by sending $q\mapsto U_1U_2$. Again $(\overline{X},X)$ is smooth. We use the embedding $\mathbb{Z}\to\mathbb{Z}^4$ which sends $n\mapsto(n,n,-n,-n)$, to define $H=\bi(\mathbb{N}^4\to(\mathbb{Z}^4/\mathbb{Z}))$. Let $h:(\mathbb{Z}^4/\mathbb{Z})\to\mathbb{Z}^2$ be the map which sends the class of $(n_1,n_2,n_3,n_4)$ to $(n_1+n_3,n_2+n_4)$, and let $K=h^{-1}(\mathbb{N}^2)$. We have an isomorphism $K\cong\mathbb{N}^2\oplus\mathbb{Z}$ by sending the class of $(n_1,n_2,n_3,n_4)$ to $(n_1+n_3,n_2+n_4,n_1-n_2)$. We thus find $$\overline{X}_1=W[H]=W[U_1,U_2,V_1,V_2]/(U_1U_2-V_1V_2),$$ $Z=W[S_1,S_2,S_3^{\pm}]$ and $g:Z\to\overline{X}_1$ is given by $U_1\mapsto S_1S_3,\quad U_2\mapsto S_2S_3^{-1},\quad V_1\mapsto S_1,\quad V_2\mapsto S_2$. Note that in this case the projections $q_j:\overline{X}_1\to\overline{X}$ are flat. Now consider the base change with $T=W[q]/q\to S$ defined by sending $q\mapsto 0$. Then, in contrast to (a), we find $\overline{X}_{1,1}=\overline{X}_{1,2}=\overline{X}_{T,1}$ (with $\overline{X}_{1,1}, \overline{X}_{1,2}, \overline{X}_{T,1}$ as in (a)).\\ 
(c) Using the criterion \ref{tglatt} one checks that the log schemes with boundary mentioned in \arabic{bspcol1}.\arabic{bspcol2}(b)--(d) and \arabic{bspbdl1}.\arabic{bspbdl2} are smooth. In fact, the example (a) just discussed is a special case of \arabic{bspcol1}.\arabic{bspcol2} (b) or \arabic{bspbdl1}.\arabic{bspbdl2}. Example (b) (or rather its base change with $T=W[q]/q\to S$ as above) is a special case of \arabic{bspcol1}.\arabic{bspcol2} (d).\\

\begin{lem} $\Omega^1_{(\overline{X},X)/T}$ is locally free of finite rank if $(\overline{X},X)$ is weakly smooth over $T$.
\end{lem}

{\sc Proof:} The same as in the classical case.\\

\section{Semistable log schemes with boundary}
\label{semist}
In this Section $k$ is a field, $Q=\mathbb{N}$ with generator $q$ and $T=(\spec(k),Q\stackrel{0}{\to}k)$.

\subsection{Definitions}

\addtocounter{satz}{1}{\bf \arabic{section}.\arabic{satz}} A {\it standard semistable $T$-log scheme with boundary} is a $T$-log scheme with boundary isomorphic to:$$(\overline{X},X)=(\spec(\frac{k[t_1,\ldots,t_{i_2}]}{(t_1,\ldots,t_{i_1})}),\spec(\frac{k[t_1,\ldots,t_{i_1},t_{i_1+1}^{\pm},\ldots,t_{i_2}^{\pm}]}{(t_1,\ldots,t_{i_1})}))$$for some integers $1\le i_1\le i_2$ such that $$P=\mathbb{N}^{i_2}\to{\cal N}_{\overline{X}}(\overline{X}),\quad 1_i\mapsto t_i\quad\mbox{ for }1\le i\le i_2$$ $$Q=\mathbb{N}\to P^{\rm gp}=\mathbb{Z}^{i_2},\quad q\mapsto (1_1,\ldots,1_{i_1},r_{i_1+1},\ldots,r_{i_2})$$with some $r_{j}\in\mathbb{Z}$ for $i_1+1\le j\le i_2$ is a chart in the sense of \arabic{defcha1}.\arabic{defcha2}. A {\it semistable $T$-log scheme with boundary} is a $T$-log scheme with boundary $(\overline{Y},Y)$ such that \'{e}tale locally on $\overline{Y}$ there exist morphisms $(\overline{Y},Y)\to(\overline{X},X)$ to standard semistable $T$-log schemes with boundary such that $\overline{Y}\to\overline{X}$ is strict and log smooth, and $Y={\overline Y}\times_{\overline{X}}X$. Note that $Y$ is then a semistable $k$-log scheme in the usual sense defined in \arabic{bspcol1}.\arabic{bspcol2}(d).

A {\it normal crossing variety over $k$} is a $k$-scheme which \'{e}tale locally admits smooth morphisms to the underlying schemes of semistable $k$-log schemes.

Following \cite{fkato} we say that a log structure ${\cal N}_{\overline{Y}}$ on a normal crossing variety $\overline{Y}$ over $k$ is {\it of embedding type} if \'{e}tale locally on $\overline{Y}$ the log scheme $(\overline{Y},{\cal N}_{\overline{Y}})$ is isomorphic to a semistable $k$-log scheme. (The point is that we do not require a {\it global} structure map of log schemes $(\overline{Y},{\cal N}_{\overline{Y}})\to T$.)\\  

\addtocounter{satz}{1}{\bf \arabic{section}.\arabic{satz}} Let us discuss for a moment the standard semistable $T$-log schemes with boundary $(\overline{X},X)$. If in the above definition $r_j\ge 0$ for all $j$, then $f:X\to T$ actually extends to a (non log smooth in general) usual morphism of log schemes $\overline{f}:\overline{X}\to T$. If even $r_j=0$ for all $j$ then $\overline{f}$ is nothing but a semistable $k$-log scheme with an additional horizontal divisor not interfering with the structure map of log structures; in particular it is log smooth. If at least $r_j\in\{0,1\}$ for all $j$ the morphism $\overline{f}$ is ideally smooth in the sense of Ogus \cite{ogid}. Examples with $r_j=1$ for all $j$ are those in \arabic{bspcol1}.\arabic{bspcol2}(d).

The concept of semistable $T$-log schemes with boundary helps us to also understand the cases with local numbers $r_j\notin\{0,1\}$: Any $(\overline{Y},Y)$ semistable $T$-log scheme with boundary is smooth, by Theorem \ref{tglatt}, and as we will see below this implies analogs of classical results for their cohomology. Examples of semistable $T$-log schemes with boundary with local numbers $r_j$ possibly not in $\{0,1\}$ are those in \arabic{bspbdl1}.\arabic{bspbdl2} or those from \ref{Grothe} below. Or think of a flat family of varieties over $\spec(k[q])$ with smooth general fibre and whose reduced subscheme of the special fibre is a normal crossing variety, but where some components of the special fibre may have multiplicities $>1$: then unions of irreducible components of this special fibre with multiplicity $=1$ are semistable $T$-log schemes with boundary. One more big class of examples with local numbers $r_j$ possibly not in $\{0,1\}$ is obtained by the following lemma, which follows from computations with local coordinates: 

\begin{lem} Let $Y\to\overline{Y}$ be an embedding of $k$-schemes which \'{e}tale locally looks like the underlying embedding of $k$-schemes of a semistable $T$-log scheme with boundary (i.e. for each geometric point $y$ of $\overline{Y}$ there is a semistable $T$-log scheme with boundary which on underlying schemes looks like $Y\to\overline{Y}$ around $y$). Suppose ${\cal N}_{\overline{Y}}$ is a log structure of embedding type on $\overline{Y}$ such that $(Y,{\cal N}_{\overline{Y}}|_Y)$ is a semistable $k$-log scheme (for an appropriate structure morphism to $T$). Then $((\overline{Y},{\cal N}_{\overline{Y}}),Y)$ is a semistable $T$-log scheme with boundary.
\end{lem} 

\addtocounter{satz}{1}{\bf \arabic{section}.\arabic{satz}} Fumiharu Kato in \cite{fkato} has worked out precise criteria for these two properties of normal crossing varieties over $k$ --- to admit a log structure of embedding type, resp. to admit a log structure of semistable type. Now suppose we are given a semistable $T$-log scheme $Y$. An ``optimal'' compactification would be a dense open embedding into a proper semistable $k$-log scheme in the classical sense, or at least into an ideally smooth proper $k$-log scheme; however, advocating the main idea of this paper, a dense open embedding $Y\to\overline{Y}$ into a log scheme $\overline{Y}$ such that $(\overline{Y},Y)$ is a proper semistable $T$-log scheme with boundary is also very useful, and this might be easier to find, or (more importantly) be naturally at hand in particular situations.\\  

\subsection{De Rham cohomology}

Here we assume $\kara(k)=0$. Let $Z$ be a smooth $k$-scheme and let $V$ be a normal crossing divisor on $Z$. Suppose there exists a flat morphism $f:(Z-V)\to\spec(k[q])$, smooth above $q\ne 0$ and with semistable fibre $X$ above the origin $q=0$. Let $\overline{X}$ be the closure of $X$ in $Z$ and suppose also that $\overline{X}\cup V$ is a normal crossing divisor on $Z$. Endow $Z$ with the log structure defined by $\overline{X}\cup V$ and endow all subschemes of $Z$ with the induced log structure (we will suppress mentioning of this log structure in our notation). Then $(\overline{X},X)$ is a semistable $T$-log scheme with boundary. Let $D=\overline{X}\cap V=\overline{X}-X$ and let $\overline{X}=\cup_{1\le i\le a}\overline{X}_i$ be the decomposition into irreducible components in a fixed ordering and suppose that each $\overline{X}_i$ is classically smooth. Let $\Omega_{X/T}^{\bullet}$ be the relative logarithmic de Rham complex.\\

\begin{pro}\label{Grothe} The restriction map$$R\Gamma(\overline{X},\Omega_{(\overline{X},X)/T}^{\bullet})\to R\Gamma(X,\Omega_{X/T}^{\bullet})$$is an isomorphism.
\end{pro}

{\sc Proof:} We use a technique of Steenbrink \cite{steen} to reduce to a standard fact. Let $\Omega_Z^{\bullet}$ be the de Rham complex over $k$ on $Z$ with logarithmic poles along $\overline{X}\cup V$. Note that $\dlog(f^*(q))\in\Gamma(Z-V,\Omega_Z^{1})$ extends uniquely to a global Section $\theta\in\Gamma(Z,\Omega_Z^{1})$. Let $\Omega_{Z,V}^{\bullet}$ be the de Rham complex on $Z$ with logarithmic poles only along $V$; thus $\Omega_{Z,V}^{\bullet}$ is a {\it subcomplex} of $\Omega_Z^{\bullet}$. Define the vertical weight filtration on $\Omega_Z^{\bullet}$ by$$P_j\Omega^i_Z=\bi(\Omega^j_Z\otimes\Omega^{i-j}_{Z,V}\to \Omega^i_Z).$$For $j\ge 1$ let $\overline{X}^j$ be the disjoint sum of all $\cap_{i\in I}\overline{X}_i$ where $I$ runs through the subsets of $\{1,\ldots,a\}$ with $j$ elements. Let $\tau_j:\overline{X}^j\to \overline{X}$ be the canonical map and let $\Omega_{\overline{X}^j}^{\bullet}$ be the de Rham complex on $\overline{X}^j$ with logarithmic poles along $\overline{X}^j\cap \tau_{j}^{-1}(D)$. Then we have isomorphisms of complexes$$(*)\quad\quad {\rm res}:\Gr_j\Omega_Z^{\bullet}\cong \tau_{j,*}\Omega_{\overline{X}^j}^{\bullet}[-j],$$characterized as follows: Let $x_1,\ldots,x_d$ be local coordinates on $Z$ such that $x_i$ for $1\le i\le a\le d$ is a local coordinate for $\overline{X}_i$. If $$\omega=\alpha\wedge\dlog(x_{i_1})\wedge\ldots \wedge\dlog(x_{i_j})\in P_j\Omega_Z^{\bullet}$$ with $i_1<\ldots<i_j<a$, then ${\rm res}$ sends the class of $\omega$ to the class of $\alpha$. Now let $$A^{pq}=\Omega_Z^{p+q+1}/P_q\Omega_Z^{p+q+1},\quad\quad\quad P_jA^{pq}=P_{2q+j+1}\Omega_Z^{p+q+1}/P_q\Omega_Z^{p+q+1}.$$Using the differentials $d':A^{pq}\to A^{p+1,q},\quad\omega\mapsto d\omega$ and $d'':A^{pq}\to A^{p,q+1},\quad\omega\mapsto \theta\wedge\omega$ we get a filtered double complex $A^{\bullet\bullet}$. We claim that $$0\to\frac{\Omega^p_Z\otimes{\cal O}_{\overline{X}}}{(\Omega^{p-1}_Z\otimes{\cal O}_{\overline{X}})\wedge\theta}\stackrel{\wedge\theta}{\to}A^{p0}\stackrel{\wedge\theta}{\to}A^{p1}\stackrel{\wedge\theta}{\to}\ldots$$is exact. Indeed, it is enough to show that for all $p$, all $j\ge 2$ the sequences$$\Gr_{j-1}\Omega^{p-1}_Z\stackrel{\wedge\theta}{\to}\Gr_{j}\Omega^{p}_Z\stackrel{\wedge\theta}{\to}\Gr_{j+1}\Omega^{p+1}_Z\stackrel{\wedge\theta}{\to}\ldots$$$$0\to P_0\Omega^{p-1}_Z/{\cal J}_{\overline{X}}.\Omega^{p-1}_Z\stackrel{\wedge\theta}{\to}\Gr_{1}\Omega^{p}_Z\stackrel{\wedge\theta}{\to}\Gr_{1}\Omega^{p+1}_Z\stackrel{\wedge\theta}{\to}\ldots$$are exact, where ${\cal J}_{\overline{X}}=\ke({\cal O}_{Z}\to{\cal O}_{\overline{X}})$. This follows from $(*)$ and the exactness of $$0\to P_0\Omega^{p}_Z/{\cal J}_{\overline{X}}.\Omega^{p}_Z\to\tau_{1,*}\Omega^p_{\overline{X}^1}\to\tau_{2,*}\Omega^p_{\overline{X}^2}\to\ldots.$$The claim follows. It implies that the maps$$\Omega_{(\overline{X},X)/T}^{p}=\frac{\Omega^p_Z\otimes{\cal O}_{\overline{X}}}{(\Omega^{p-1}_Z\otimes{\cal O}_{\overline{X}})\wedge\theta}\to A^{p0}\subset A^p,\quad\omega\mapsto (-1)^p\theta\wedge\omega$$define a quasi-isomorphism $\Omega_{(\overline{X},X)/T}^{\bullet}\to A^{\bullet}$, hence a spectral sequence$$E_1^{-r,q+r}=H^q(\overline{X},\Gr_r A^{\bullet})\Longrightarrow H^q(\overline{X},\Omega_{(\overline{X},X)/T}^{\bullet}).$$Now we can of course repeat all this on $Z-V$ instead of $Z$, and restriction from $Z$ to $Z-V$ gives a canonical morphism between the respective spectral sequences. That this is an isomorphism can be checked on the initial terms, and using the isomorphism $(*)$ this boils down to proving that the restriction maps$$H^p(\overline{X}^j,\Omega_{\overline{X}^j}^{\bullet})\to H^p(X^j,\Omega_{\overline{X}^j}^{\bullet})$$are isomorphisms where we set $X^j=\overline{X}^j\cap\tau_j^{-1}(X)$. But this is well known. The proof is finished.\\

\addtocounter{satz}{1}{\bf \arabic{section}.\arabic{satz}} Now assume ${\overline X}$ is proper. Similar to the classical Hodge theory, the Hodge filtration on  $$H^p(\overline{X},\Omega_{(\overline{X},X)/T}^{\bullet})=H^p(X,\Omega_{X/T}^{\bullet})$$ obtained by stupidly filtering $\Omega_{(\overline{X},X)/T}^{\bullet}$ should be meaningful. Another application of Proposition \ref{Grothe} might be a Poincar\'{e} duality theorem. Suppose the underlying scheme of $\overline{X}$ is of pure dimension $d$. Let ${\cal I}_{D}=\ke({\cal O}_{\overline{X}}\to{\cal O}_{D})$ and define the de Rham cohomology with compact support of $(\overline{X},X)/T$ as $$R\Gamma(\overline{X},{\cal I}_D\otimes\Omega_{(\overline{X},X)/T}^{\bullet}).$$It is a natural question to ask if this is dual to $R\Gamma(\overline{X},\Omega_{(\overline{X},X)/T}^{\bullet})=R\Gamma(X,\Omega_{X/T}^{\bullet})$. The key would be as usual the construction of a trace map $H^d(\overline{X},{\cal I}_D\otimes\Omega_{(\overline{X},X)/T}^d)\to k$.\\

\addtocounter{satz}{1}{\bf \arabic{section}.\arabic{satz}} Another application of semistable $T$-log schemes with boundary is the possibility to define the notion of {\it regular singularities} of a given integrable log connection on a semistable $T$-log scheme $X$, provided we have an embedding $X\to\overline{X}$ such that $(\overline{X},X)$ is a proper semistable $T$-log scheme with boundary.\\

\addtocounter{satz}{1}{\bf \arabic{section}.\arabic{satz}} Here is an application of the construction in \arabic{bspbdl1}.\arabic{bspbdl2} to the de Rham cohomology of certain semistable $k$-log schemes (a simplified variant of the application given in \cite{hkstrat}; in fact, the present paper formalizes and generalizes a key construction from \cite{hkstrat}). In \arabic{bspbdl1}.\arabic{bspbdl2} assume that $\kara(k)=0$ and that $M$ is the intersection of {\it all} irreducible components of $Y$. Recall that we constructed a morphism of log schemes $V_M\to S=(\spec(k[q]),1\mapsto q)$. For $k$-valued points $\alpha\to S$ (with pull back log structure) let $V_{M}^{\alpha}=V_M\times_S{\alpha}$. Using the $S$-log scheme with boundary $(P_M,V_M)$ one can show that the derived category objects $R\Gamma(V_{M}^{\alpha},\Omega^{\bullet}_{V_M^{\alpha}/\alpha})$ (with $\Omega^{\bullet}_{V_M^{\alpha}/\alpha}$ the relative logarithmic de Rham complex; if $\alpha\ne 0$ this is the classical one) are canonically isomorphic for varying $\alpha$. Namely, the canonical restriction maps $$R\Gamma(P_{M},\Omega^{\bullet}_{(P_M,V_M)/S})\to R\Gamma(V_{M}^{\alpha},\Omega^{\bullet}_{V_M^{\alpha}/\alpha})$$are isomorphisms for all $\alpha$.\\

\subsection{Crystalline cohomology}

Let $\tilde{S}$ be a scheme such that ${\cal O}_{\tilde{S}}$ is killed by a non-zero integer, $I\subset{\cal O}_{\tilde{S}}$ a quasi-coherent ideal with DP-structure $\gamma$ on it, and let $\tilde{\cal L}$ be a fine log structure on $\tilde{S}$. Let $(S,{\cal L})$ be an exact closed log subscheme of $(\tilde{S},\tilde{\cal L})$ defined by a sub-DP-ideal of $I$ and let $f:(X,{\cal N})\to(S,{\cal L})$ be a log smooth and integral morphism of log schemes. An important reason why log crystalline cohomology of $(X,{\cal N})$ over $(\tilde{S},\tilde{\cal L})$ works well is that locally on $X$ there exist smooth and integral, hence flat morphisms $\tilde{f}:(\tilde{X},\tilde{\cal N})\to(\tilde{S},\tilde{\cal L})$ with $f=\tilde{f}\times_{(\tilde{S},\tilde{\cal L})}(S,{\cal L})$. This implies that the crystalline complex of $X/\tilde{S}$ (with respect to any embedding system) is flat over ${\cal O}_{\tilde{S}}$, see \cite{hyoka} 2.22, and on this property many fundamental theorems rely.

Now let $W$ be a discrete valuation ring of mixed characteristic $(0,p)$ with maximal ideal generated by $p$. For $n\in \mathbb{N}$ let $W_n=W/(p^n)$, $k=W_1$ and $K_0=\quot(W)$, and let $T_n$ be the exact closed log subscheme of $S=W[Q]$ defined by the ideal $(p^n,q)$ (abusing previous notation we now take $\spec(W)$ as the base scheme $W$ of \arabic{basnot1}.\arabic{basnot2}). Thus $T=T_1$. We will often view $T$-log schemes with boundary as $T_n$-log schemes with boundary for $n\in \mathbb{N}$.\\

\begin{lem}\label{loclif} Let $(\overline{Y},Y)/T$ be a semistable $T$-log scheme with boundary. Then there exist \'{e}tale locally on $\overline{Y}$ smooth $T_n$-log schemes with boundary $(\overline{Y}_n,Y_n)$ such that $(\overline{Y},Y)=(\overline{Y}_n,Y_n)\overline{\times}_{T_n}T$, the closed immersion $(\overline{Y},Y)\to(\overline{Y}_n,Y_n)$ is boundary exact, and such that $\Omega^1_{(\overline{Y}_n,Y_n)/T_n}$ is flat over ${\cal O}_{T_n}$ and commutes with base changes $T_m\to T_n$ for $m\le n$.
\end{lem}

{\sc Proof:} We may suppose that there is a strict and log smooth morphism $$h:(\overline{Y},Y)\to (\overline{X},X)=(\spec(\frac{k[t_1,\ldots,t_{i_2}]}{(t_1,\ldots,t_{i_1})}),\spec(\frac{k[t_1,\ldots,t_{i_1},t_{i_1+1}^{\pm},\ldots,t_{i_2}^{\pm}]}{(t_1,\ldots,t_{i_1})}))$$for some integers $1\le i_1\le i_2$ such that $P=\mathbb{N}^{i_2}$ is a chart for $(\overline{X},X)$ sending $1_i\mapsto t_i$ for $1\le i\le i_2$ and such that the structure map is given by $$Q=\mathbb{N}\to P^{\rm gp}=\mathbb{Z}^{i_2},\quad q\mapsto (1_1,\ldots,1_{i_1},r_{i_1+1},\ldots,r_{i_2})$$with some $r_{j}\in\mathbb{Z}$ for $i_1+1\le j\le i_2$. We lift $(\overline{X},X)$ to$$(\overline{X}_n,X_n)=(\spec(\frac{W_n[t_1,\ldots,t_{i_2}]}{(t_1,\ldots,t_{i_1})}),\spec(\frac{W_n[t_1,\ldots,t_{i_1},t_{i_1+1}^{\pm},\ldots,t_{i_2}^{\pm}]}{(t_1,\ldots,t_{i_1})}))$$using the same formulas for the log structure maps. Local liftings of $h$ to $(\overline{X}_n,X_n)$ result from the classical theory, since `strict and log smooth' is equivalent to `smooth on underlying schemes'.\\

\begin{lem}\label{kriskey} Let $n\in\mathbb{N}$ and let $(\overline{Y},Y)\to(\overline{X}_i,X_i)$ be boundary exact closed immersions into smooth $T_n$-log schemes with boundary ($i=1,2$). Then there exist \'{e}tale locally on $(\overline{X}_1\overline{\times}_{T_n}\overline{X}_2)$ factorizations of the diagonal embedding $$(\overline{Y},Y)\stackrel{\iota}{\to}(\overline{Z},Z)\stackrel{}{\to}(\overline{X}_1\overline{\times}_{T_n}\overline{X}_2,X_1\times_{T_n}X_2)$$with $\iota$ a boundary exact closed immmersion, the map $\overline{Z}\to\overline{X}_1\overline{\times}_{T_n}\overline{X}_2$ log \'{e}tale, the projections $p_i:\overline{Z}\to\overline{X}_i$ strict and log smooth, and with the following property: Let $\overline{D}_{12}$ (resp. $\overline{D}_i$) denote the ${\rm DP}$ envelopes of (the underlying scheme morphism of) $\overline{Y}\to\overline{Z}$ (resp. of $\overline{Y}\to\overline{X}_i$), and let $q_i:\overline{D}_{12}\to\overline{D}_i$ be the canonical projections. Then there exist $u_{i1},\ldots,u_{im_i}\in {\cal O}_{\overline{D}_{12}}$ for $i=1$ and $i=2$ such that $du_{i1},\ldots,du_{im_i}$ form a basis of $\Omega^1_{\overline{Z}/\overline{X}_i}$ and such that the assignments $U_{ij}^{[k]}\mapsto u_{ij}^{[k]}$ ($k\in\mathbb{N}$) induce isomorphisms$$q_i^{-1}{\cal O}_{\overline{D}_i}\langle U_{i1},\ldots,U_{im_i}\rangle\cong{\cal O}_{\overline{D}_{12}}$$where on the left hand side we mean the ${\rm DP}$ envelope of the free polynomial ring.
\end{lem} 

\addtocounter{satz}{1}{\bf \arabic{section}.\arabic{satz}}\newcounter{defcry1}\newcounter{defcry2}\setcounter{defcry1}{\value{section}}\setcounter{defcry2}{\value{satz}} Lemma \ref{kriskey} follows from Proposition \ref{keykoh}, and the same proofs give variants of Proposition \ref{keykoh} and Lemma \ref{kriskey} for more than two embeddings $(\overline{Y},Y)\to(\overline{X}_i,X_i)$ (and hence with products with more than two factors). As in \cite{kalo} one shows that the DP envelopes of $(\overline{Y},Y)$ in chosen exactifications of these products (e.g. the DP envelope $\overline{D}_{12}$ in Lemma \ref{kriskey}) are independent of the chosen  exactifications. For a given semistable $T$-log scheme with boundary $(\overline{Y},Y)$ we now define its crystalline cohomology relative to $T_n$ by the standard method (cf. \cite{hyoka} 2.18): Choose an open covering $\overline{Y}=\cup_{\overline{U}\in\overline{\cal U}}\overline{U}$ and for each $(\overline{U},U=Y\cap\overline{U})$ a lift $(\overline{U}_n,U_n)$ as in Lemma \ref{loclif}. Taking products we get a simplicial $T_n$-log scheme with boundary $(\overline{U}^{\bullet}_n,U^{\bullet}_n)$ which is an embedding system for $(\overline{Y},Y)$ over $T_n$. Let $\overline{D}_n^{\bullet}$ be the DP envelope of $(\overline{Y},Y)$ in $(\overline{U}_n^{\bullet},U^{\bullet}_n)$, i.e. the simplicial scheme formed by the DP envelopes of local exactifications of $(\overline{Y},Y)\to(\overline{U}_n^{\bullet},U^{\bullet}_n)$ as in Lemma \ref{kriskey}. Then we set$$R\Gamma_{\rm crys}((\overline{Y},Y)/T_n)=R\Gamma(\overline{D}_n^{\bullet},\Omega^{\bullet}_{(\overline{U}_n,U_n)/T_n}\otimes{\cal O}_{\overline{D}_n^{\bullet}}).$$That this definition is independent of the chosen embedding follows from Lemma \ref{kriskey} and the DP Poincar\'{e} lemma.\\

\begin{lem} (a) For $m\le n$ we have $$R\Gamma_{\rm crys}((\overline{Y},Y)/T_m)\cong R\Gamma_{\rm crys}((\overline{Y},Y)/T_n)\otimes^{\mathbb{L}}_{W_n}W_m.$$(b) If $\overline{Y}$ is proper over $k$, the cohomology of $$R\lim_{\stackrel{\leftarrow}{n}}R\Gamma_{\rm crys}((\overline{Y},Y)/T_n)$$ (resp. of $R\Gamma_{\rm crys}((\overline{Y},Y)/T_n)$) is finitely generated over $W$ (resp. over $W_n$). 
\end{lem}

{\sc Proof:} Just as in \cite{hyoka} 2.22 one deduces from Lemmata \ref{loclif} and \ref{kriskey} that $\Omega^{\bullet}_{(\overline{U}_n,U_n)/T_n}\otimes{\cal O}_{\overline{D}_n^{\bullet}}$ is a $W_n$-flat sheaf complex on $\overline{D}_n^{\bullet}$ and this implies (a). If $\overline{Y}$ is proper over $k$ it follows that $R\Gamma_{\rm crys}((\overline{Y},Y)/T_1)=R\Gamma(\overline{Y},\Omega^{\bullet}_{(\overline{Y},Y)/T_1})$ has finite dimensional cohomology over $k$ since each $\Omega^{j}_{(\overline{Y},Y)/T_1}$ is coherent. Together with (a) we conclude as in the classical case.\\

\addtocounter{satz}{1}{\bf \arabic{section}.\arabic{satz}} Ogus \cite{ogcon} and Shiho \cite{shiho} have defined logarithmic convergent cohomology in great generality and ``in crystalline spirit''. Here we content ourselves with the following definition. Let $E$ be a fine $T$-log scheme. Let $T_{\infty}$ be the formal log scheme $(\spf(W),1\mapsto 0)$. Choose an exact closed immersion $E\to G$ into a log smooth formal $T_{\infty}$-log scheme $G$ topologically of finite type over $W$. Associated to $G$ is a $K_0$-rigid space $G_{K_0}$ together with a specialization map ${\rm sp}$ to the special fibre of $G$. The preimage ${\rm sp}^{-1}(E)=]E[_G$ of the embedded $E$, the tube of $E$, is an admissible open subspace of $G_{K_0}$. The logarithmic de Rham complex $\Omega^{\bullet}_{G/T_{\infty}}$ on $G$ gives rise, tensored with $\mathbb{Q}$, to a sheaf complex $\Omega^{\bullet}_{G_{K_0}/T_{\infty,K_0}}$ on $G_{K_0}$ and we set$$R\Gamma_{\rm conv}(E/T_{\infty})=R\Gamma(]E[_G,\Omega^{\bullet}_{G_{K_0}/T_{\infty,K_0}}),$$an object in the derived category of $K_0$-vector spaces. If there are embeddings $E\to G$ as above only locally on $E$, one works with embedding systems.\\

Now let $Y$ be a semistable $k$-log scheme with smooth irreducible components and let $M$ be the intersection of some of its irreducible components. Endow $M$ with the structure of $T$-log scheme induced from $Y$. Note that $M$ is not log smooth over $T$ (unless $Y$ has only a single irreducible component) and its usual log crystalline cohomology is pathological; it does not provide a canonical integral lattice in the log convergent cohomology of $M$, as we will now construct one by another method. In \arabic{bspbdl1}.\arabic{bspbdl2} we constructed a $S_1$-log scheme with boundary $(P_M,V_M)$ where $S_1$ is the exact closed log subscheme of $S$ defined by the ideal $(p)$. Perform the base change with the exact closed subscheme $T$ of $S_1$ defined by the ideal $(q)$ to get $(P_M^0,V_M^0)=(P_M\times_{S_1}T,V_M\times_{S_1}T)$. This is a semistable $T$-log scheme with boundary as defined above.\\

\begin{satz}\label{crisconv} There exists a canonical isomorphism $$R\lim_{\stackrel{\leftarrow}{n}}R\Gamma_{\rm crys}(({P_M^0},V_M^0)/T_n)\otimes_{W}K_0\cong R\Gamma_{\rm conv}(M/T_{\infty}).$$In particular, if $M$ is proper, each $R^j\Gamma_{\rm conv}(M/T_{\infty})$ is finite dimensional.\\ 
\end{satz}

{\sc Proof:} {\it Step 1:} The map is$$R\lim_{\stackrel{\leftarrow}{n}}R\Gamma_{\rm crys}(({P_M^0},V_M^0)/T_n)\otimes_{W}K_0\to R\lim_{\stackrel{\leftarrow}{n}}R\Gamma_{\rm crys}(({V_M^0},V_M^0)/T_n)\otimes_{W}K_0$$$$=R\lim_{\stackrel{\leftarrow}{n}}R\Gamma_{\rm crys}(V_M^0/T_n)\otimes_{W}K_0\stackrel{(i)}{\cong}R\Gamma_{\rm conv}(V_M^0/T_{\infty})\to R\Gamma_{\rm conv}(M/T_{\infty})$$where the left hand side in $(i)$ is the usual log crystalline cohomology of $V_M^0/T_n$ and the isomorphism $(i)$ holds by log smoothness of $V_M^0/T$. That this map is an ismorphism can be checked locally.\\{\it Step 2:} We may therefore assume that there exists a smooth (in the classical sense) affine connected $\spec(W)$-scheme $\tilde{\cal M}=\spec(\tilde{B})$ lifting $M$ and that the invertible sheaves ${\cal F}_j|_M$ on $M$ are trivial (notation from \arabic{bspbdl1}.\arabic{bspbdl2}); let $v_j$ be a generator of ${\cal F}_j|_M$. Furthermore we may assume that the divisor $D$ on $M$ (the intersection of $M$ with all irreducible components of $Y$ not containing $M$) lifts to a (relative $\spec(W)$) normal crosssings divisor $\tilde{\cal D}$ on $\tilde{\cal M}$. Let $$\tilde{\cal V}_{\cal M}=\spec(\tilde{B}[x_j]_{j\in I})$$$$\tilde{\cal P}_{\cal M}=\times_{\tilde{\cal M}}(\proj(\tilde{B}[y_j,x_j]_{j\in I}).$$Identifying the free variable $x_j$ with a lift of $v_j$ we view $\tilde{\cal V}_{\cal M}$ as a lift of $V_M$; identifying moreover the free variable $y_j$ with a lift of $1_{{\cal O}_M}$ we view $\tilde{\cal P}_{\cal M}$ as a lift of $P_M$; identifying a homogenous element $s\in\tilde{B}[x_j]_{j\in I}$ of degree $n$ with the degree zero element $s/y_j^n$ of $\tilde{B}[y_j^{\pm},x_j]$ we view $\tilde{\cal V}_{\cal M}$ as an open subscheme of $\tilde{\cal P}_{\cal M}$. As in \arabic{bspbdl1}.\arabic{bspbdl2} we factor the distinguished element $a\in\sym_{{\cal O}_M}(\oplus({\cal F}_j)_{j\in I})(M)$ as $a=a_0.(\oplus_{j\in I}v_j)$ with defining equation $a_0\in{\cal O}_M$ of the divisor $D$ in $M$. Lift $a_0$ to a defining equation $\tilde{a}_0\in \tilde{B}$ of $\tilde{\cal D}$ in $\tilde{\cal M}$. This $\tilde{a}_0$ also defines a normal crossing divisor  $\tilde{\cal D}_{\tilde{\cal V}}$ on $\tilde{\cal V}_{\cal M}$. Set $\tilde{a}=\tilde{a}_0\prod_{j\in I}x_j\in \tilde{B}[x_j]_{j\in I}$ and consider the following normal crossing divisor on $\tilde{\cal P}_{\cal M}$: the union of $\tilde{\cal P}_{\cal M}-\tilde{\cal V}_{\cal M}$ with the closure (in $\tilde{\cal P}_{\cal M}$) of the zero set of $\tilde{a}$ (in $\tilde{\cal V}_{\cal M}$). It defines a log structure on $\tilde{\cal P}_{\cal M}$. Define a morphism $\tilde{\cal V}_{\cal M}\to S$ by sending $q\mapsto \tilde{a}$. We have constructed a lift of the $S_1$-log scheme with boundary $(P_M,V_M)$ to a $S$-log scheme with boundary $(\tilde{\cal P}_{\cal M},\tilde{\cal V}_{\cal M})$. Moreover, if we denote by $\tilde{\cal T}_{\infty}$ the exact closed log subscheme of $S$ defined by the ideal $(q)$, then the $\tilde{\cal T}_{\infty}$-log scheme with boundary $(\tilde{\cal P}_{\cal M}^0,\tilde{\cal V}_{\cal M}^0)=(\tilde{\cal P}_{\cal M}\times_S\tilde{\cal T}_{\infty},\tilde{\cal V}_{\cal M}\times_S\tilde{\cal T}_{\infty})$ is a lift of the $T$-log scheme with boundary $(P_M^0,V_M^0)$.\\{\it Step 3:} Denote by ${\cal P}_{\cal M}^0$ (resp. ${\cal V}_{\cal M}^0$, resp. ${\cal M}$, resp. ${\cal D}_{\cal V}$) the $p$-adic completions of $\tilde{\cal P}_{\cal M}^0$ (resp. of $\tilde{\cal V}_{\cal M}^0$, resp. of $\tilde{\cal M}$, resp. of $\tilde{\cal D}_{\tilde{\cal V}}$). Denote by ${\cal P}_{{\cal M},n}^0$ (resp. ${\cal V}_{{\cal M},n}^0$, resp. ${\cal M}_n$, resp. ${\cal D}_{{\cal V},n}$) the reduction modulo $p^n$. Let $\Omega^{\bullet}_{{\cal P}_{\cal M}^0/T_{\infty}}$ be the $p$-adic completion of the de Rham complex of the $\tilde{\cal T}_{\infty}$-log scheme with boundary $(\tilde{\cal P}_{\cal M}^0,\tilde{\cal V}_{\cal M}^0)$. Its reduction $\Omega^{\bullet}_{{\cal P}_{\cal M}^0/T_{\infty}}\otimes(\mathbb{Z}/p^n)$ modulo $p^n$ is the de Rham complex $\Omega^{\bullet}_{{\cal P}_{{\cal M},n}^0/T_n}$ of the ${T}_{n}$-log scheme with boundary $({\cal P}_{{\cal M},n}^0,\tilde{\cal V}_{{\cal M},n}^0)$. Observe that the differentials on $\Omega^{\bullet}_{{\cal P}_{\cal M}^0/T_{\infty}}$ pass to differentials on $\Omega^{\bullet}_{{\cal P}_{\cal M}^0/T_{\infty}}\otimes_{{\cal O}_{{\cal P}_{\cal M}^0}}{\cal O}_{\cal M}$ where we use the zero section ${\cal M}\to{\cal P}_{\cal M}^0$. Let $\Omega^{\bullet}_{{\cal P}_{\cal M}^0/T_{\infty}}\otimes\mathbb{Q}$ be the complex on the rigid space ${\cal P}_{{\cal M},K_0}^0$ obtained by tensoring with $K_0$ the sections of $\Omega^{\bullet}_{{\cal P}_{\cal M}^0/T_{\infty}}$ over open affine pieces of ${\cal P}_{\cal M}^0$. Similarly define $\Omega^{\bullet}_{{\cal P}_{\cal M}^0/T_{\infty}}\otimes_{{\cal O}_{{\cal P}_{\cal M}^0}}{\cal O}_{\cal M}\otimes\mathbb{Q}$. By definition we have $$R\Gamma_{\rm conv}(M/T_{\infty})=R\Gamma(]M[_{{\cal P}_{\cal M}},\Omega^{\bullet}_{{\cal P}_{\cal M}^0/T_{\infty}}\otimes\mathbb{Q}),$$
$$R\lim_{\stackrel{\leftarrow}{n}}R\Gamma_{\rm crys}(({P_M^0},V_M^0)/T_n)\otimes_{W}K_0=R\lim_{\stackrel{\leftarrow}{n}}R\Gamma({\cal P}_{{\cal M},n}^0,\Omega^{\bullet}_{{\cal P}_{{\cal M},n}^0/T_n})\otimes_{W}K_0.$$In view of$$R\Gamma(]M[_{{\cal P}_{\cal M}},\Omega^{\bullet}_{{\cal P}_{\cal M}^0/T_{\infty}}\otimes_{{\cal O}_{{\cal P}_{\cal M}^0}}{\cal O}_{\cal M}\otimes\mathbb{Q})=R\lim_{\stackrel{\leftarrow}{n}}R\Gamma({\cal P}_{{\cal M},n}^0,\Omega^{\bullet}_{{\cal P}_{{\cal M},n}^0/T_n}\otimes_{{\cal O}_{{\cal P}_{{\cal M},n}^0}}{\cal O}_{{\cal M}_n})\otimes_{W}K_0$$
it is therefore enough to show that the maps$$f_n:R\Gamma({\cal P}_{{\cal M},n}^0,\Omega^{\bullet}_{{\cal P}_{{\cal M},n}^0/T_n})\to R\Gamma({\cal P}_{{\cal M},n}^0,\Omega^{\bullet}_{{\cal P}_{{\cal M},n}^0/T_n}\otimes_{{\cal O}_{{\cal P}_{{\cal M},n}^0}}{\cal O}_{{\cal M}_n})$$
$$g:R\Gamma(]M[_{{\cal P}_{\cal M}},\Omega^{\bullet}_{{\cal P}_{\cal M}^0/T_{\infty}}\otimes\mathbb{Q})\to R\Gamma(]M[_{{\cal P}_{\cal M}},\Omega^{\bullet}_{{\cal P}_{\cal M}^0/T_{\infty}}\otimes_{{\cal O}_{{\cal P}_{\cal M}^0}}{\cal O}_{\cal M}\otimes\mathbb{Q})$$are isomorphisms.\\{\it Step 4:} Let ${\cal D}_{{\cal V},n}=\cup_{l\in L}{\cal D}_{n,l}$ be the decomposition of ${\cal D}_{{\cal V},n}$ into irreducible components. Let ${\cal E}'_{n}$ be the closed subscheme of ${\cal V}_{{\cal M},n}^0$ defined by $\prod_{j\in I}x_j\in \Gamma({\cal V}_{{\cal M},n}^0,{\cal O}_{{\cal V}_{{\cal M},n}^0})$ and let ${\cal E}_{n}$ be the closure of ${\cal E}'_{n}$ in ${\cal P}_{{\cal M},n}^0$. Let ${\cal E}_{n}=\cup_{j\in I}{\cal E}_{n,j}$ be its decomposition into irreducible components. For a pair $P=(P_I,P_L)$ of subsets $P_I\subset I$ and $P_L\subset L$ let$${\cal G}_P=(\cap_{j\in P_I}{\cal E}_{n,j})\cap(\cap_{l\in P_L}{\cal D}_{n,l}),$$so we drop reference to $n$ in our notation, for convenience. Also for convenience we denote the sheaf complex $\Omega^{\bullet}_{{\cal P}_{{\cal M},n}^0/T_n}$ on ${\cal P}_{{\cal M},n}^0$ simply by $\Omega^{\bullet}$. For two pairs $P, P'$ as above with $P_I\cup P_L\ne \emptyset$, with $P_I\subset P'_I$ and $P_L=P'_{L}$ consider the canonical map $$w_{P,P'}:\Omega^{\bullet}\otimes{\cal O}_{{\cal G}_P}\to\Omega^{\bullet}\otimes{\cal O}_{{\cal G}_{P'}}$$of sheaf complexes on ${\cal P}_{{\cal M},n}^0$. We claim that the map $R\Gamma({\cal P}_{{\cal M},n}^0,w_{P,P'})$ induced by $w_{P,P'}$ in cohomology is an isomorphism. For this we may of course even assume $P_I'=P_I\cup\{j_0\}$ for some $j_0\in I$, $j_0\notin P_I$. In the ${\cal O}_{{\cal G}_{P'}}$-module $\Omega^1\otimes{\cal O}_{{\cal G}_{P'}}$ we fix a complement $N$ of the submodule generated by (the class of) $\dlog(x_{j_0})\in\Gamma({\cal P}_{{\cal M},n}^0,\Omega^1\otimes{\cal O}_{{\cal G}_{P'}})$ as follows. We use the identification$$\frac{(\Omega^1_{\tilde{\cal M}}(\log(\tilde{\cal D}))\otimes{\cal O}_{{\cal G}_{P'}})\oplus(\oplus_{j\in I}{\cal O}_{{\cal G}_{P'}}.\dlog(x_j))}{{\cal O}_{{\cal G}_{P'}}.\dlog(\tilde{a})}=\Omega^1\otimes{\cal O}_{{\cal G}_{P'}}$$(with $\Omega^1_{\tilde{\cal M}}(\log(\tilde{\cal D}))$ the differential module of $(\tilde{\cal M},(\mbox{log str. def. by }\tilde{\cal D}))\to(\spec(W),\mbox{triv.})$). If $P_L\ne\emptyset$ we may assume that we can factor our $\tilde{a}_0\in\tilde{B}$ from above as $\tilde{a}_0=\tilde{a}_0'h$ with $h\in\tilde{B}$ whose zero set in $\tilde{\cal M}=\spec(\tilde{B})$ reduces modulo $p^n$ to an irreducible component of $\cup_{l\in P_L}{\cal D}_{l,n}$. We may assume that the ${\cal O}_{\tilde{\cal M}}$-submodule of $\Omega^1_{\tilde{\cal M}}(\log(\tilde{\cal D}))$ generated by $\dlog(h)$ admits a complement $N'$. Then we get the isomorphism $$(N'\otimes{\cal O}_{{\cal G}_{P'}})\oplus(\oplus_{j\in I}{\cal O}_{{\cal G}_{P'}}.\dlog(x_j))\cong \Omega^1\otimes{\cal O}_{{\cal G}_{P'}}$$(use $\dlog(\tilde{a})=\dlog(h)+\dlog({\tilde{a}}/{h})$). If there exists $j'\in P_I$ we get the isomorphism$$(\Omega^1_{\tilde{\cal M}}(\log(\tilde{\cal D}))\otimes{\cal O}_{{\cal G}_{P'}})\oplus(\oplus_{j\in I-\{j'\}}{\cal O}_{{\cal G}_{P'}}.\dlog(x_j))\cong\Omega^1\otimes{\cal O}_{{\cal G}_{P'}}$$(use $\dlog(\tilde{a})=\dlog(x_{j'})+\dlog({\tilde{a}}/{x_{j'}})$). In both cases, dropping the $j_0$-summand in the left hand side we get $N$ as desired. We see that the ${\cal O}_{{\cal G}_{P'}}$-subalgebra $N^{\bullet}$ of $\Omega^{\bullet}\otimes{\cal O}_{{\cal G}_{P'}}$ generated by $N$ is stable for the differential $d$, and that we have $$\Omega^{\bullet}\otimes{\cal O}_{{\cal G}_{P'}}=N^{\bullet}\otimes_{W_n}C^{\bullet}$$
as complexes, where $C^{\bullet}$ is the complex $C^0=W_n$, $C^1=W_n.\dlog(x_{j_0})$ (here $\dlog(x_{j_0})$ is nothing but a symbol), $C^m=0$ for $m\ne 0,1$, and zero differential. Let ${\cal R}=\proj (W_n[y_{j_0},x_{j_0}])$. We have a canonical map ${\cal G}_P\to{\cal R}$. Let $D^{\bullet}$ be the ${\cal O}_{\cal R}$-subalgebra of $\Omega^{\bullet}\otimes{\cal O}_{{\cal G}_{P}}$ generated by $\dlog(x_{j_0})\in\Gamma({\cal P}_{{\cal M},n}^0,\Omega^1\otimes{\cal O}_{{\cal G}_{P}})$. It is stable for the differential $d$, and we find$$\Omega^{\bullet}\otimes{\cal O}_{{\cal G}_{P}}=N^{\bullet}\otimes_{W_n}D^{\bullet}$$as complexes, where $N^{\bullet}$ is mapped to $\Omega^{\bullet}\otimes{\cal O}_{{\cal G}_{P}}$ via the natural map (and section of $w_{P,P'}$)$$\Omega^{\bullet}\otimes{\cal O}_{{\cal G}_{P'}}\to\Omega^{\bullet}\otimes{\cal O}_{{\cal G}_{P}}$$ induced by the structure map ${\cal E}_{n,j_0}\to{\cal M}_n$. This map also induces a map $C^{\bullet}\to D^{\bullet}$, and it is enough to show that the latter induces isomorphisms in cohomology. But $$H^m({\cal P}_{{\cal M},n}^0,D^{\bullet})\cong H^m({\bf P}^1_{W_n},\Omega_{{\bf P}^1_{W_n}}^{\bullet}(\log\{0,\infty\})),$$which is $W_n$ if $0\le m\le 1$ and zero otherwise, because of the degeneration of the Hodge spectral sequence (\cite{kalo} 4.12) and $\Omega^1_{{\bf P}^1}(\log\{0,\infty\})\cong{\cal O}_{{\bf P}^1}$. So $C^{\bullet}$ and $D^{\bullet}$ have the same cohomology.\\{\it Step 5:} We now show that $f_n$ is an isomorphism. Let ${\cal F}_I=\cup_{j\in I}{\cal E}_{n,j}$, let ${\cal F}_L=\cup_{l\in L}{\cal D}_{n,l}={\cal D}_{{\cal V},n}$ and ${\cal F}_{I,L}={\cal F}_I\cap{\cal F}_L$. All the following tensor products are taken over ${\cal O}_{{\cal P}_{{\cal M},n}^0}$. We will show that in$$\Omega^{\bullet}=\Omega^{\bullet}\otimes {\cal O}_{{\cal P}_{{\cal M},n}^0}\stackrel{\alpha}{\longrightarrow}\Omega^{\bullet}\otimes {\cal O}_{{\cal F}_I}\stackrel{\beta}{\longrightarrow}\Omega^{\bullet}\otimes {\cal O}_{{\cal G}_{(I,\emptyset)}}=\Omega^{\bullet}\otimes{\cal O}_{{\cal M}_n}$$both $\alpha$ and $\beta$ induce isomorphisms in cohomology. The exact sequences
$$0\longrightarrow {\cal O}_{{\cal P}_{{\cal M},n}^0}\longrightarrow {\cal O}_{{\cal F}_J}\oplus {\cal O}_{{\cal F}_L}\longrightarrow {\cal O}_{{\cal F}_{I,L}}\longrightarrow0$$$$0\longrightarrow{\cal O}_{{\cal F}_{I}}\longrightarrow{\cal O}_{{\cal F}_{I}}\oplus  {\cal O}_{{\cal F}_{I,L}} \longrightarrow {\cal O}_{{\cal F}_{I,L}}\longrightarrow0$$ show that, to prove that $\alpha$ induces cohomology isomorphisms, it is enough to prove that $\Omega^{\bullet}\otimes{\cal O}_{{\cal F}_L}\to\Omega^{\bullet}\otimes {\cal O}_{{\cal F}_{I,L}}$ induces cohomology isomorphisms. To see this, it is enough to show that both $\Omega^{\bullet}\otimes {\cal O}_{{\cal F}_L}\stackrel{\gamma}{\to}\Omega^{\bullet}\otimes {\cal O}_{{\cal F}_L\cap {\cal G}_{(I,\emptyset)}}$ and $\Omega^{\bullet}\otimes {\cal O}_{{\cal F}_{I,L}}\stackrel{\delta}{\to}\Omega^{\bullet}\otimes {\cal O}_{{\cal F}_L\cap {\cal G}_{(I,\emptyset)}}$ induce cohomology isomorphisms. Consider the exact sequence\begin{gather}0\longrightarrow {\cal O}_{{\cal F}_L}\longrightarrow\oplus_{l\in L}{\cal O}_{{\cal G}_{(\emptyset,\{l\})}}\longrightarrow\oplus_{\stackrel{L'\subset L}{|L'|=2}}{\cal O}_{{\cal G}_{(\emptyset,L')}}\longrightarrow\ldots\longrightarrow{\cal O}_{{\cal G}_{(\emptyset,L)}} \longrightarrow 0\tag{$*$}\end{gather}Comparison of the exact sequences $(*)\otimes\Omega^{\bullet}$ and $(*)\otimes {\cal O}_{{\cal F}_L\cap {\cal G}_{(I,\emptyset)}}\otimes\Omega^{\bullet}$ shows that to prove that $\gamma$ induces cohomology isomorphisms, it is enough to show this for $\Omega^{\bullet}\otimes{\cal O}_{{\cal G}_{(\emptyset,L')}}\to\Omega^{\bullet}\otimes{\cal O}_{{\cal G}_{(I,L')}}$ for all $\emptyset\ne L'\subset L$; but this has been done in Step 2. Comparison of $(*)\otimes {\cal O}_{{\cal F}_{I,L}}\otimes\Omega^{\bullet}$ and $(*)\otimes {\cal O}_{{\cal F}_L\cap {\cal G}_{(I,\emptyset)}}\otimes\Omega^{\bullet}$ shows that to prove that $\delta$ induces cohomology isomorphisms, it is enough to show this for $\Omega^{\bullet}\otimes{\cal O}_{{\cal F}_I\cap{\cal G}_{(\emptyset,L')}}\stackrel{\epsilon_G}{\to}\Omega^{\bullet}\otimes {\cal O}_{{\cal G}_{(I,L')}}$ for all $\emptyset\ne L'\subset L$. Consider the exact sequence \begin{gather}0\longrightarrow {\cal O}_{{\cal F}_I}\longrightarrow\oplus_{j\in I}{\cal O}_{{\cal G}_{(\{j\},\emptyset)}}\longrightarrow\oplus_{\stackrel{I'\subset I}{|I'|=2}}{\cal O}_{{\cal G}_{(I',\emptyset)}}\longrightarrow\ldots\longrightarrow{\cal O}_{{\cal G}_{(I,\emptyset)}} \longrightarrow 0\tag{$**$}\end{gather} The exact sequence $(**)\otimes {\cal O}_{{\cal F}_I\cap{\cal G}_{(\emptyset,L')}}\otimes\Omega^{\bullet}$ shows that to prove that $\epsilon_G$ induces cohomology isomorphisms, it is enough to show this for $\Omega^{\bullet}\otimes {\cal O}_{{\cal G}_{(I',L')}}\to\Omega^{\bullet}\otimes {\cal O}_{{\cal G}_{(I,L')}}$ for all $\emptyset\ne I'\subset {I}$; but this has been done in Step 2. The exact sequence $(**)\otimes\Omega^{\bullet}$ shows that to prove that $\beta$ induces cohomology isomorphisms, it is enough to show this for $\Omega^{\bullet}\otimes {\cal O}_{{\cal G}_{(I',\emptyset)}}\to\Omega^{\bullet}\otimes {\cal O}_{{\cal G}_{(I,\emptyset)}}$
for all $\emptyset\ne I'\subset {I}$; but this has been done in Step 2. The proof that $f_n$ is an isomorphism is complete.

The proof that $g$ is an isomorphism is essentially the same: While Step 4 above boiled down to $H^m({\bf P}^1_{W_n},\Omega_{{\bf P}^1_{W_n}}^{\bullet}(\log\{0,\infty\}))=W_n$ if $0\le m\le 1$, and $=0$ for other $m$, one now uses $H^m({\bf D}^0_{K_0},\Omega_{{\bf D}^0_{K_0}}^{\bullet}(\log\{0\}))=K_0$ if $0\le m\le 1$, and $=0$ for other $m$ (here ${\bf D}^0_{K_0}$ is the open unit disk over $K_0$). The formal reasoning from Step 5 is then the same. The theorem is proven.\\

\addtocounter{satz}{1}{\bf \arabic{section}.\arabic{satz}}\newcounter{uniconv1}\newcounter{uniconv2}\setcounter{uniconv1}{\value{section}}\setcounter{uniconv2}{\value{satz}} Also unions $H$ of irreducible components of $Y$ are not log smooth over $T$ (unless $H=Y$) and their usual log crystalline cohomology is not useful. However, if $H^{\heartsuit}$ denotes the complement in $H$ of the intersection of $H$ with the closure of $Y-H$ in $Y$, then $(H,H^{\heartsuit})$ is a semistable $T$-log scheme with boundary. There is natural map$$h:R\Gamma_{\rm conv}(H/T)\longrightarrow R\lim_{\stackrel{\leftarrow}{n}}R\Gamma_{\rm crys}((H,H^{\heartsuit})/T_n)\otimes_{W}K_0,$$constructed as follows. We say a $T_{\infty}$-log scheme is strictly semistable if all its irreducible components are smooth $W$-schemes and if \'{e}tale locally it is the central fibre of a morphism $\spec(W[t_1,\ldots,t_{n}])\to \spec(W[t]),\quad t\mapsto t_1\cdots t_{m}$ (some $n\ge m\ge1$), with the log structures defined by the vanishing locus of $t$ resp. of $t_1\cdots t_{m}$. We find an \'{e}tale cover $Y=\{Y_i\}_{i\in I}$ of $Y$ and for each $i\in I$ a semistable $T_{\infty}$-log scheme ${\mathcal Y}_i$ together with an isomorphism ${\mathcal Y}_i\otimes_Wk\cong Y_i$. Taking suitable blowing ups in the products of these ${\mathcal Y}_i$ (a standard procedure, compare for example \cite{mokr}) we get an embedding system for $Y$ over $T_{\infty}$ where a typical local piece $Y_J=\prod_Y(Y_i)_{i\in J}$ of $Y$ is exactly embedded as $Y_J\to {\mathcal Y}_J$ with ${\mathcal Y}_J$ a semistable $T_{\infty}$-log scheme and such there is a closed subscheme ${\mathcal H}_J$ of ${\mathcal Y}_J$, the union of some of its irreducible components, such that $H\times_Y Y_J={\mathcal H}_J\times_{{\mathcal Y}_J}Y$. Now ${\mathcal Y}_J$ is log smooth over $T_{\infty}$, hence its $p$-adic completion $\widehat{\mathcal Y}_J$ may be used to compute $R\Gamma_{\rm conv}(H\times_Y Y_J/T)$. On the other hand, let ${\mathcal H}_J^{\heartsuit}\subset {\mathcal H}_J$ be the open subscheme which is the complemet in ${\mathcal Y}_J$ of all irreducible components of ${\mathcal Y}_J$ which are not fully contained in ${\mathcal H}_J$. Then $({\mathcal H}_J,{\mathcal H}_J^{\heartsuit})$ is a smooth $T_{\infty}$-log scheme with boundary, hence its $p$-adic completion may be used to compute $R\lim_{\stackrel{\leftarrow}{n}}R\Gamma_{\rm crys}((H\times_Y Y_J,H^{\heartsuit}\times_Y Y_J)/T_n)\otimes_{W}K_0$. By the proof of \cite{berfi} Proposition 1.9 there is a natural map from the structure sheaf of the tube $]H\times_Y Y_J[_{\widehat{\mathcal Y}_J}$ to the structure sheaf of the $p$-adically completed DP envelope, tensored with $\mathbb{Q}$, of $H\times_Y Y_J$ in ${\mathcal H}_J$. It induces a map between our de Rham complexes in question, hence we get $h$. By the same local argument which showed the isomorphy of the map $g$ in the proof of Theorem \ref{crisconv} we see that $h$ is an isomorphism; the work on local lifts of $P_M^0$ there is replaced by work on local lifts of $Y$ here. In particular, if $H$ is proper, each $R^j\Gamma_{\rm conv}(H/T_{\infty})$ is finite dimensional.\\

\addtocounter{satz}{1}{\bf \arabic{section}.\arabic{satz}} Suppose $k$ is perfect. Then there is a canonical Frobenius endomorphism on the log scheme $T_n$ (cf. \cite{hyoka} 3.1): The canonical lift of the $p$-power map on $k$ to an endomorphism of $W_n$, together with the endomorphism of the log structure which on the standard chart $\mathbb{N}$ is multiplication with $p$. We can also define a Frobenius endomorphism on $R\Gamma_{\rm crys}((\overline{Y},Y)/T_n)$ for a semistable $T$-log scheme with boundary $(\overline{Y},Y)$, because we can define a Frobenius endomorphism on the embedding system used in \arabic{defcry1}.\arabic{defcry2}, compatible with that on $T_n$. Namely, on a standard $T_n$-log scheme with boundary $(\overline{X}_n,X_n)$ as occurs in the proof of Lemma \ref{loclif} we act on the underlying scheme by the Frobenius on $W_n$ and by $t_i\mapsto t_i^p$ (all $i$), and on the log structure we act by the unique compatible map which on our standard chart $\mathbb{N}^{i_2}$ is multiplication with $p$. Then we lift these endomorphisms further (using the lifting property of classical smoothness) to Frobenius lifts of our $\overline{Y}$-covering and hence to the embedding system.\\

\addtocounter{satz}{1}{\bf \arabic{section}.\arabic{satz}} We finish with perspectives on possible further developments.\\(1) Mokrane \cite{maghreb} defines the crystalline cohomology of a classically smooth $k$-scheme $U$ as the log crystalline cohomology with poles in $D$ of a smooth compactification $X$ of $U$ with $D=X-U$ a normal crossing divisor. This is a cohomology theory with the usual good properties (finitely generated, Poincar\'{e} duality, mixed if $k$ is finite). He shows that under assumptions on resolutions of singularities, this cohomology theory indeed only depends on $U$. We suggest a similar approach to define the crystalline cohomology of a semistable $k$-log scheme $U$: Compactify it (if possible) into a proper semistable $T$-log scheme with boundary $(X,U)$ and take the crystalline cohomology of $(X,U)$.

Similarly, classical rigid cohomology as defined by Berthelot \cite{berco} works with compactifications. Also here, to define log versions it might be useful to work with log schemes with boundary to avoid hypotheses on existence of compactifications by genuine log morphisms.\\ (2) We restricted our treatment of crystalline cohomology to that of semistable $T$-log schemes with boundary $(\overline{Y},Y)$ relative to $T_n$. For deformations of $T=T_1$ other then $T_n$ --- for example, $(\spec(W_n), 1\mapsto p)$ --- we have at present no suitable analogs of Lemma \ref{loclif}. However, such analogs also seem to lack in idealized log geometry: for an ideally log smooth $T$-log scheme (like the union of some irreducible components of a semistable $k$-log scheme in the usual sense), there seems to be in general no lift to a flat and ideally log smooth $(\spec(W_n), 1\mapsto p)$-log scheme. Some more foundational concepts need to be found.

Let us nevertheless propose some tentative definitions of crystalline cohomology for more general fine log schemes $T$ and more general $T$-log schemes with boundary (without claiming any results). Suppose that $p$ is nilpotent in ${\cal O}_W$ and let $(I,\delta)$ be a quasicoherent DP ideal in ${\cal O}_W$. All DP structures on ideals in ${\cal O}_W$-algebras are required to be compatible with $\delta$. Let $T_0$ be a closed subscheme of $T$ and let $\gamma$ be a DP structure on the ideal of $T_0$ in $T$. Let $(\overline{X},{X})$ be a ${T}$-log scheme with boundary, and let $\overline{X}_0$ be the closure in $\overline{X}$ of its locally closed subscheme ${X}\times_{T}{T}_0$. We say $\gamma$ extends to $(\overline{X},{X})$ if there is a DP structure $\alpha$ on the ideal of $\overline{X}_0$ in $\overline{X}$, such that the structure map ${X}\to{T}$ is a DP morphism (if $\alpha$ exists, it is unique, because ${\cal O}_{\overline{X}}\to i_*{\cal O}_{X}$ is injective). Then we say $(\overline{X},{X})$ is a $\gamma$-${T}$-log scheme with boundary. For a $\gamma$-${T}$-log scheme $(\overline{X},{X})$ we can define the crystalline site and the crystalline cohomology of $(\overline{X},{X})$ over $T$ as in the case of usual log schemes.

{\it Example:} Let $\underline{T}_0\subset \underline{T}$ be a closed immersion. Suppose ${T}$ is the DP envelope of $\underline{T}_0$ in $\underline{T}$ and ${T}_0\subset {T}$ is the closed subscheme defined by its DP ideal; we have $T_0=\underline{T}_0$ if $\delta$ extends to $\underline{T}$. Now if $(\underline{\overline{X}},\underline{X})$ is a $\underline{T}$-log scheme with boundary, we obtain a $\gamma$-${T}$-log scheme with boundary $(\overline{X},X)$ by taking as $\overline{X}$ the DP envelope of the schematic closure of the subscheme $\underline{X}\times_{\underline{T}}\underline{T}_0$ of $\overline{\underline{X}}$.\\


\begin{flushleft}
\textsc{Mathematisches Institut der Universit\"at M\"unster\\ Einsteinstrasse 62, 48149 M\"unster, Germany}\\
\textit{E-mail address}: klonne@math.uni-muenster.de
\end{flushleft}
\end{document}